\documentclass[11pt]{amsart}

\usepackage{amssymb}
\usepackage{graphicx}
\usepackage[scanall]{psfrag}

\newcommand{\bpf}{
\begin{proof}
}
\newcommand{\epf}{
\end{proof}
}

\newtheorem{thm}{Theorem}[section]
\newtheorem{cor}[thm]{Corollary}
\newtheorem{lem}[thm]{Lemma}
\newtheorem{prop}[thm]{Proposition}

\newtheorem{case}{Case}

\newtheorem{claim*}{Claim}

\theoremstyle{definition}

\newcommand{\ve}{\varepsilon}

\newcommand{\ben}{\begin{enumerate}}
\newcommand{\een}{\end{enumerate}}
\newcommand{\ita}{\item[(a)]}
\newcommand{\itb}{\item[(b)]}
\newcommand{\itc}{\item[(c)]}

\newcommand{\wh}{\widehat}

\newcommand{\mc}{\mathcal}
\newcommand{\intr}{\text{int}\,}

\newcommand{\bg}{\overline{G}}

\newcommand{\bgt}{\overline{G}_T}
\newcommand{\bgs}{\overline{G}_S}
\newcommand{\bgpp}{\overline{G}_P}

\newcommand{\ove}{\overline{e}}

\newcommand{\id}{\text{id}}

\newcommand{\dgt}{\deg_{{G}_T}}
\newcommand{\dgs}{\deg_{{G}_S}}
\newcommand{\dbgt}{\deg_{\overline{G}_T}}
\newcommand{\dbgs}{\deg_{\overline{G}_S}}

\newcommand{\Figw}[4]{
\includegraphics[width=#1]{#2}
\caption{ #3 \label{#4} } }

\begin{document}

\title[Toroidal and Klein bottle boundary slopes]
{Toroidal and Klein bottle boundary slopes}%

\author[L. G. Valdez-S\'anchez]{Luis G. Valdez-S\'anchez}
\address{Department of Mathematical Sciences,
University of Texas at El Paso\\
El Paso, TX 79968, USA}
\email{valdez@math.utep.edu}%


\date{\today. {\bf Preliminary version, update v3}}

\begin{abstract}
Let $M$ be a compact, connected, orientable, irreducible 3-manifold
and $T_0$ an incompressible torus boundary component of $M$ such
that the pair $(M,T_0)$ is not cabled. By a result of C.\ Gordon, if
$(S,\partial S),(T,\partial T)\subset(M,T_0)$ are incompressible
punctured tori with boundary slopes at distance
$\Delta=\Delta(\partial S,\partial T)$, then $\Delta\leq 8$, and the
cases where $\Delta=6,7,8$ are very few and classified. We give a
simplified proof of this result (or rather, of its {\it reduction
process}), using an improved estimate for the maximum possible
number of mutually parallel negative edges in the graphs of
intersection of $S$ and $T$. We also extend Gordon's result by
allowing either $S$ or $T$ to be an essential Klein bottle.
\end{abstract}

\maketitle


\section{Introduction}\label{intro}

Let $M$ be a compact, connected, orientable, irreducible 3-manifold,
and $T_0$ an incompressible torus boundary component of $M$. If
$r_1,r_2$ are two slopes in $T_0$, we denote their {\it distance},
ie their minimum geometric intersection number in $T_0$, by
$\Delta(r_1,r_2)$. By a {\it surface} we mean a compact
2-dimensional manifold, not necessarily orientable. A properly
embedded surface in $M$ with nonempty boundary which is not a disk
is said to be {\it essential} if it is geometrically incompressible
and boundary incompressible in $M$. We will use the notion of a {\it
cabled pair $(M,T_0)$} in the sense of
\cite{gordonlith}.

Let $(F,\partial F)\subset (M,T_0)$ be a punctured torus. We say
that {\it $F$ is generated by a (an essential) Klein bottle} if
there is a (an essential, resp.) punctured Klein bottle $(P,\partial
P)\subset (M,T_0)$ such that $F$ is isotopic in $M$ to the frontier
of a regular neighborhood of $P$ in $M$. We also say that $F$ is
{\it $\mc{K}$-incompressible} if $F$ is either incompressible or
generated by an essential Klein bottle. In this paper we give a
proof of the following result.

\begin{thm}\label{main}
Let $(F_1,\partial F_1),(F_2,\partial F_2)\subset (M,T_0)$ be
$\mc{K}$-incompressible tori, and let $\Delta=\Delta(\partial
F_1,\partial F_2)$. If the pair $(M,T_0)$ is not cabled then
$\Delta\leq 8$, and if $\Delta\geq 6$ then $|\partial F_1|,|\partial
F_2|\leq 2$.
\end{thm}

The corollary below follows immediately from Theorem~\ref{main};
along with \cite[Theorem 1.2 and \S 6]{valdez6}, it can be used to
obtain the classification of the manifolds $M$ that contain
essential punctured Klein bottles with boundary slopes at distance
$\Delta\geq 6$.

\begin{cor}\label{maincor}
Let $(F_1,\partial F_1),(F_2,\partial F_2)\subset (M,T_0)$ be
punctured essential Klein bottles, and let $\Delta=\Delta(\partial
F_1,\partial F_2)$. If the pair $(M,T_0)$ is not cabled then
$\Delta\leq 8$, and if $\Delta\geq 6$ then $|\partial
F_1|=1=|\partial F_2|$, with $\Delta=6,8$.\hfill\qed
\end{cor}

Theorem~\ref{main} is well known when the surfaces $F_{\alpha}$ are
both tori, in which case it follows from the proof of
\cite[Proposition 1.5]{gordon5}. The case where both surfaces are
Klein bottles has been discussed more recently in \cite[Corollary
1.5]{lee1} (for $\Delta\geq 5$) and
\cite[Theorem 1.4]{mati1} (for $\Delta\geq 5$), under the added
hypothesis that $M$ is hyperbolic. Thus, for $\Delta\geq 6$, modulo
the classification of the manifolds $M$, Theorem~\ref{main} and its
corollary extend the range of applicability of
\cite[Proposition 1.5]{gordon5} to include the case of
essential Klein bottles, and of
\cite[Corollary 1.5]{lee1} and \cite[Theorem 1.4]{mati1} to allow
for manifolds that may not be hyperbolic.

A general approach to the proof of results similar to
Theorem~\ref{main} involves what we may call a {\it reduction
process,} where, say, a condition on the distance between the
boundary slopes, like $\Delta\geq 6$, creates `large' families of
parallel edges, whose presence may restrict the number of boundary
components of at least one surface to be `small', or the topology of
$M$ to be `degenerate', in some sense. If the `small' cases are
sufficiently small, they can be dealt with separately or classified
completely. In fact, for $\Delta\geq 6$, combining the
classification of the pairs $(M,T_0)$ in
\cite[Proposition 1.5]{gordon5} with Theorem~\ref{main} and
\cite[Theorem 1.2]{valdez6}, it follows that there are exactly
four manifolds $(M,T_0)$ in Theorem~\ref{main}, all obtained via
Dehn fillings along one of the boundary components of the Whitehead
link exterior, and that if $\Delta=6,8$ and $F_\alpha$ is a torus
then $F_{\alpha}$ is incompressible and generated by a once
punctured Klein bottle.

In the proof of Theorem~\ref{main} we present here we use some
fundamental results from the paper \cite{gordonlith}, with the
addition of Lemma 2.1 \cite[\S 2]{gordon5} (on parallelism of
edges), the notion of jumping number \cite[\S 2]{gordon5}, and the
parity rules from \cite{cgls,tera1,tera10}; the new ingredients are
contained in Proposition~\ref{prop1}, the main technical result of
this paper, which roughly states that if $(M,T_0)$ is not cabled and
contains two $\mc{K}$-incompressible tori $(T,\partial T),
(T',\partial T')\subset (M,T_0)$ with $\Delta(\partial T,\partial
T')\geq 1$, then, for any surface $S\subset M$ that intersects $T$
in essential graphs, any collection of mutually parallel negative
edges of the graph $S\cap T\subset S$ has at most $|\partial T|+1$
edges, unless $M$ is one of three exceptional toroidal manifolds, in
which case $\Delta(\partial T,\partial T')=1,2$ or 4. We remark that
the current best bound used in similar contexts is $2\cdot|\partial
T|$, for $t\geq 4$ (cf \cite[Corollary 5.5]{gordon5}). It is the use
of the upper bound $|\partial T|+1$ of Proposition~\ref{prop1} that
gives rise to a rather short reduction process for
Theorem~\ref{main}.

The paper is organized as follows. In Section~\ref{secA} we present
several basic definitions and facts related to the graphs of
intersection produced by two surfaces in $M$ with transverse
intersection. Section~\ref{neg} is devoted to the discussion of
bounds for the sizes of collections of mutually parallel edges in
the graphs of intersection of two surfaces in $M$; the first two
subsections deal with the case of positive edges and some known
facts for the case of negative edges, and the remaining two sections
contain the proof of Proposition~\ref{prop1}. Finally, the proof of
Theorem~\ref{main} is given in Section~\ref{delta6}.

We thank Masakazu Teragaito and Sangyop Lee for their careful
reading of preliminary versions of this preprint and many helpful
suggestions.

\section{Preliminaries}\label{secA}

Let $M$ be a compact, connected, orientable, irreducible 3-manifold
with an incompressible torus boundary component $T_0$. For any
nontrivial slope $r\subset T_0$, $M(r)$ will denote the Dehn filled
manifold $M\cup_{T_0} V$, where $V$ is a solid torus such that $r$
bounds a disk in $V$. If $F\subset M$ is a properly embedded surface
and $r$ is the slope of the circles $F\cap T_0$, then $\wh{F}$ will
denote the surface in $M(r)$ obtained from $F$ by capping off any
components of $\partial F$ in $T_0$ with disjoint meridian disks in
$V$.

Let $F_1,F_2$ be any two properly embedded surfaces in $M$
(orientable or not) which intersect transversely in a minimum number
of components; in particular, if $r_{\alpha}$ is the slope of the
circles $\partial F_{\alpha}\cap T_0$ in $T_0$, and
$\Delta=\Delta(r_1,r_2)$, then any two components of $\partial
F_1\cap T_0$ and $\partial F_2\cap T_0$ intersect transversely in
$\Delta$ points.

We say that $G_{F_1}=F_1\cap F_2\subset F_1$ and $G_{F_2}=F_1\cap
F_2\subset F_2$ are the {\it graphs of intersection} between $F_1$
and $F_2$. Either of these graphs is {\it essential} if each
component of $F_1\cap F_2$ is geometrically essential in the
corresponding surface. The graph $G_{F_{\alpha}}$ has {\it fat
vertices} the components of $\partial F_{\alpha}$ and edges the arc
components of $F_1\cap F_2$; there may also be some circle
components present. An edge of $F_1\cap F_2$ with both endpoints in
$T_0$ is called an {\it internal edge}.

Let $n_1=|\partial F_1\cap T_0|$ and $n_2=|\partial F_2\cap T_0|$.
We label the components of $\partial F_{\alpha}\cap T_0$ as
$\partial_1F_{\alpha},\partial_2F_{\alpha},\dots,
\partial_{n_{\alpha}}F_{\alpha}$,
consecutively in their order of appearance along $T_0$ (in some
direction), and then label each intersection point between
$\partial_i F_1$ and $\partial_j F_2$ with $j$ in $G_{F_1}$ and $i$
in $G_{F_2}$. In this way, any endpoint of an edge of $F_1\cap F_2$
that lies in $T_0$ gets a label in each graph of intersection, and
internal edges get labels at both endpoints.

Following \cite{tera1,tera10}, we orient the components of $\partial
F_{\alpha}\cap T_0$ coherently on $T_0$, and say that an internal
edge $e$ of $F_1\cap F_2$ has a {\it positive} or {\it negative
sign} in $G_{F_{\alpha}}$ depending on whether the orientations of
the components of $\partial F_{\alpha}$ (possibly the same) around a
small rectangular regular neighborhood of $e$ in $F_{\alpha}$ appear
as in Fig.~\ref{n20}.
\begin{figure}
\psfrag{e}{$e$}
\psfrag{positive edge}{Positive edge}
\psfrag{negative edge}{Negative edge}
\Figw{2.5in}{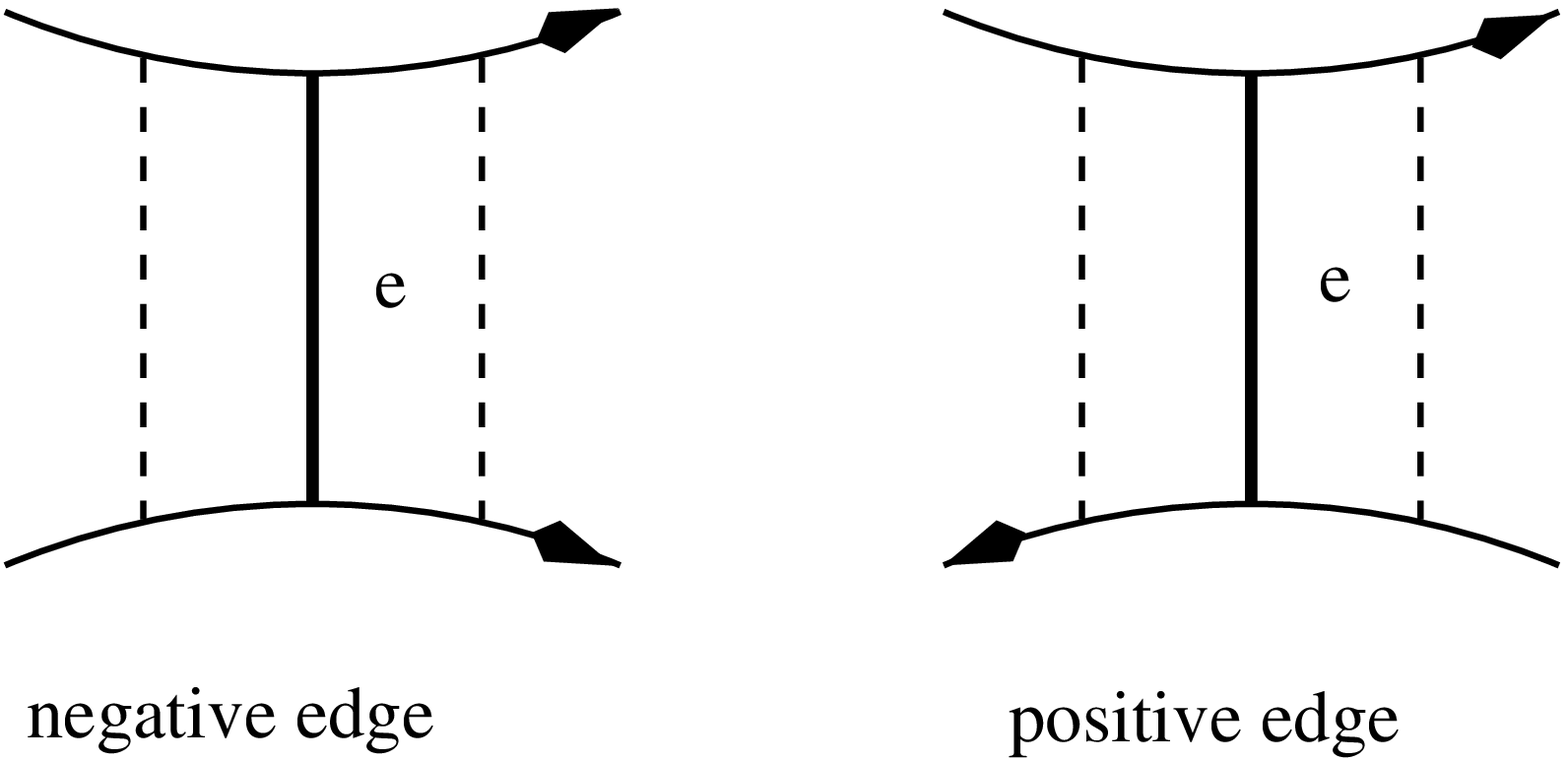}{}{n20}
\end{figure}

Alternatively (cf \cite{cgls}), if $F_{\alpha}$ is orientable, we
fix an orientation on $F_{\alpha}$, induce an orientation on the
components of $\partial F_{\alpha}\cap T_0$, and then say that two
components of $\partial F_{\alpha}\cap T_0$ have {\it the same
parity} if their given orientations agree on $T_0$, and {\it
opposite parity} otherwise. This divides the components of $\partial
F_{\alpha}\cap T_0$ into two parity classes, and we may call the
vertices in one class {\it positive}, and the vertices in the other
class {\it negative}. It is then not hard to see that an internal
edge of $F_1\cap F_2$ is positive (negative) in $F_{\alpha}$ iff it
connects two vertices of $G_{F_{\alpha}}$ of the same (opposite,
resp.) parity. In this context, if the vertices
$\partial_iF_{\alpha}$ of $F_{\alpha}$ are all of the same parity we
will say that $F_{\alpha}$ is {\it polarized}, and that it is {\it
neutral} if there are the same number of vertices of either parity.

A collection of edges in $G_{F_{\alpha}}$  whose union is a circle
in $\wh{F}_{\alpha}$ (where the circle is constructed in the obvious
way, by collapsing the vertices of $F_{\alpha}$ on $T_0$ into points
in $\wh{F}_{\alpha}$) is called a {\it cycle}. A cycle in
$F_{\alpha}$ is {\it nontrivial} if it is not contained in a disk in
$\wh{F}_{\alpha}$. We call a cycle in $F_{\alpha}$ consisting of a
single edge a {\it loop edge\/;} notice that if $F_{\alpha}$ is
orientable then a loop edge in $F_{\alpha}$ is positive.

Two edges of $F_1\cap F_2$ are said to be {\it parallel} in
$F_{\alpha}$ if they cobound a rectangular disk subregion in
$F_{\alpha}$. Suppose that two internal edges $e,e'$ of $F_1\cap
F_2$ are positive, parallel, and consecutive in $F_{\alpha}$, and
let $F$ be the disk face in $G_{F_{\alpha}}$ they cobound. We say
that $F$ is an {\it S-cycle face of type $\{j,j+1\}$} of
$G_{F_{\alpha}}$ (with $j,j+1$ well defined mod $n_{\beta}$) if the
labels at the endpoints of each edge $e,e'$ are $j$ and $j+1$; this
is a restricted version of the more general notion of a Scharlemann
cycle, which we will not use in this paper.

The following lemma summarizes several fundamental results we will
use in the sequel.

\begin{lem}\label{mainref}
Let $F_1,F_2$ be properly embedded surfaces in $M$ with essential
graphs of intersection $G_{F_{1}},  G_{F_{2}}$.

\ben
\ita
Parity Rule (\cite{cgls,tera1,tera10}): for
$\{\alpha,\beta\}=\{1,2\}$, an internal edge of $F_1\cap F_2$ is
positive in $G_{F_{\alpha}}$ iff it is negative in $G_{F_{\beta}}$.

\itb
Suppose $(M,T_0)$ is not cabled and $F_1,F_2$ are orientable. Then
no two internal edges of $F_1\cap F_2$ are parallel in both
$G_{F_1}$ and $G_{F_2}$ (\cite[Lemma 2.5]{gordon5}), and if
$n_{\alpha}\geq 2$ and $E$ is a family of mutually parallel,
consecutive, internal negative edges in $G_{F_{\beta}}$ then no
component of $F_{\alpha}\setminus\cup E$ is a disk in $F_{\alpha}$
(\cite{gordonlith}).

\itc
If $(M,T_0)$ is not cabled, $F_1$ is planar, $F_2$ is toroidal, and
$\partial F_1,\partial F_2\subset T_0$, then $\Delta\leq 5$
(\cite{gordonlith}).
\hfill\qed
\een
\end{lem}

\subsection{Reduced graphs}

Let $G$ be an essential connected graph on a compact punctured
surface $\mc{F}$, of the type constructed above. We let $V(G),E(G)$
denote the sets of (fat) vertices and edges of $G$, respectively.
Cutting each edge of $G$ along some interior point splits the edges
into pieces which we call the {\it local edges of $G$}. The degree
of a vertex $v$ of $G$, denoted by $\deg_G(v)$ or $\deg(v)$, is then
the number of local edges of $G$ that are incident to $v$.

For an integer $k\geq 0$, the notation $\deg\geq k$ ($\deg\equiv k$)
in $G$ will mean that $\deg(v)\geq k$ ($\deg(v)= k$, resp.) holds
for any $v\in V(G)$. Thus, the degree of any vertex $\partial_i
F_{\alpha}$ of $G_{F_{\alpha}}$ is $\Delta\cdot n_{\beta}$ and the
labels $1,2,\dots,n_{\beta}$ repeat $\Delta$ times in blocks
consecutively around $\partial_iF_{\alpha}$.

Let $N(E(G))$ be a small product neighborhood of $E(G)$ in $\mc{F}$.
Then the closure of any component of $\mc{F}\setminus N(E(G))$ is
called a {\it face} of $G$. Observe that if $F$ is any face of $G$,
then $\partial F$ is a union of segments of the form $e\times
0,e\times 1$'s, called the {\it edges of $F$}, and segments coming
from the $\partial_i \mc{F}$'s, called the {\it corners of $F$}. We
call a disk face of $G$ with $n$ sides (and $n$ corners) a {\it disk
$n$-face;} disk $2$-faces or $3$-faces are also referred to as {\it
bigons or triangles}, respectively.

The graph $G$ is said to be {\it reduced} if no two of its edges are
parallel. The {\it reduced graph $\bg$ of $G$} is the graph obtained
by amalgamating any maximal collection of mutually parallel edges of
$G$ into a single edge. Notice that any disk face in a reduced graph
is at least a triangle.

The next result gives two useful facts about reduced graphs on a
torus.

\begin{lem}\label{3v}
Let $G$ be a reduced graph on a torus with $V$ vertices, $E$ edges,
and $\deg\geq 1$.
\ben
\ita
If $\deg\geq 6$ in $G$ then $\deg\equiv 6$ in $G$ and all faces of
$G$ are triangles.

\itb
If $G$ has no triangle faces then $G$ has a vertex of degree at most
$4$.
\een
\end{lem}

\bpf
Part (a) is well known (cf \cite[Lemma 3.2]{gordon5}). For part (b),
let $d$ be the number of disk faces of $G$ and set $n=\min\{\deg(u)
\ | \ u\text{ is a vertex of } G\}\geq 1$. Then $nV\leq 2E$, and
since any disk face of $G$ is at least a $4$-face then $4d\leq 2E$.
Combining these relations with Euler's relation $E\leq V+d$ then
implies that $n\leq 4$, hence $G$ has a vertex of degree at most 4.
\epf

\subsection{Edge orbits and permutations}\label{orbits}

We will denote any edge in the reduced graph $\bg_{F_{\alpha}}$
generically by the symbol $\ove$. Hence, $\ove$ represents a
collection $e_1,e_2,\dots,e_k$ of mutually parallel, consecutive,
same sign edges in $G_{F_{\alpha}}$, in which case we say that
$|\ove|=k$ is the {\it size of $\ove$}, and that the sign of $\ove$
is positive (negative) if all the edges in $\ove$ are positive
(negative, resp).

Suppose that $n_{\beta}\geq 2$, and that $E$ is a collection of
$n_{\beta}$ mutually parallel, consecutive internal edges of
$G_{F_{\alpha}}$. We assume that these edges have endpoints in the
vertices $u_i,u_{i'}$ of $G_{F_{\alpha}}$ (with $u_i=u_{i'}$
allowed), and that all edges in $\ove$ are oriented to run from
$u_i$ to $u_{i'}$ (the orientation is arbitrary if $u_i=u_{i'}$).
Then each of the labels $1,2\dots,n_{\beta}$ appears exactly once at
the endpoints of the edges of $E$ at each of the vertices $u_i$ and
$u_{i'}$, and so the set $E$ {\it induces} a permutation $\sigma$ on
the set $\{1,2\dots,n_{\beta}\}$, defined by matching the labels at
the endpoints of the edges of $E$ in $u_i$ with the corresponding
labels at the endpoints of these edges in $u_{i'}$. This permutation
is of the form $\sigma(x)\equiv\alpha-\ve\cdot x\mod n_{\beta}$,
where $\ve=+1,-1$ is the sign of the edges in $E$ (see
Figs.~\ref{n03}(a) and \ref{n01}); reversing the orientation of the
edges replaces $\sigma$ with its inverse. Observe that if the edges
in $E$ are positive then $\sigma^2=\id$, and that $\sigma\neq\id$
whenever $F_{\beta}$ is orientable by the parity rule.

More generally, it is not hard to see that if $E'$ is any collection
of mutually parallel, consecutive internal edges of
$G_{F_{\alpha}}$, with $|E'|\geq n_{\beta}$, then any two
subfamilies of $E'$ with $n_{\beta}$ consecutive edges induce the
same permutation; we refer to this common permutation as the
permutation {\it induced} by $E'$.

The union in $G_{F_{\beta}}$ of all edges in $E$, along with all
vertices of $G_{F_{\beta}}$ at their endpoints, form a subgraph
$\Gamma_E$ of $G_{F_{\beta}}$; we call any component of $\Gamma_E$
an {\it edge orbit of $E$}. Each orbit of $\sigma$ then corresponds
uniquely to some edge orbit of $E$: for the labels of the vertices
of $G_{F_{\beta}}$ at the endpoints of the edges in an edge orbit of
$E$ form an orbit of $\sigma$.

\subsection{Strings}

We denote by $I_{i,i+1}$ the annulus cobounded in $T_0$ by the
circles $\partial_i F_{\alpha},\partial_{i+1} F_{\alpha}$, with
labels $i,i+1$ well defined mod $n_{\beta}$, and call it a {\it
string of $F_{\alpha}$}.

Notice that the corners of any face of $G_{F_{\beta}}$ are spanning
arcs along some of the strings of $F_{\alpha}$. For $F_{\alpha}$ an
orientable surface, let $N(F_{\alpha})=F_{\alpha}\times[0,1]$ be a
small product regular neighborhood of $F_{\alpha}$ in $M$; if $F$ is
a face of $G_{F_{\beta}}$, we will say that $F$ {\it locally lies on
one side of $F_{\alpha}$} if $F$ intersects only one of the two
surfaces $F_{\alpha}\times 0$ or $F_{\alpha}\times 1$.

\subsection{$\mc{K}$-incompressible tori}\label{ktori}

Suppose that the punctured torus $(T,\partial T)\subset (M,T_0)$ is
generated by an essential punctured Klein bottle $P\subset M$, and
that $S\subset M$ is a properly embedded surface which intersects
$P$ in essential graphs $G_{S,P}=S\cap P\subset S$ and $G_P=S\cap
P\subset P$. Let $N(P)$ be a regular neighborhood of $P$ in $M$, and
isotope $T$ so that $T=\text{fr}\, N(P)$. For $N(P)$ small enough,
the intersection $S\cap T$ will be transverse and the graphs
$G_{S,T}=S\cap T\subset S$ and $G_T=S\cap T\subset T$ will also be
essential; in fact, the graph $G_{S,T}$ will be the frontier in $S$
of the regular neighborhood $N(P)\cap S$ of all the components of
$G_{S,P}$. Moreover, if $\ove$ is an edge of $\bgpp$, then $\ove$
gives rise to two distinct edges $\ove_1,\ove_2$ in $\bgt$, each of
the same size as $\ove$, which are parallel in $N(P)$, and if
$\partial S\subset T_0$ and $|\ove|\geq |\partial S|$, then the
edges $\ove,\ove_1$, and $\ove_2$ all have the same sign and induce
the same permutation.

In particular, if $T_1,T_2$ are $\mc{K}$-incompressible tori in
$(M,T_0)$, then it is possible to isotope $T_1$ or $T_2$ so that
both graphs of intersection $G_{T_1}$ and $G_{T_2}$ are essential.

\subsection{S-cycles and Klein bottles}

In this section we assume that $(T,\partial T)\subset (M,T_0)$ is a
twice punctured torus and $S$ is a properly embedded surface in $M$
which intersects $T$ in essential graphs $G_S,G_T$. In particular,
all edges of $S\cap T$ are internal, and if $G_S$ has an S-cycle
face then $T$ is neutral by the parity rule.

The next result follows in part from the proof of \cite[Lemma
5.2]{gordonlu4}; we include a sketch of its proof for the
convenience of the reader.

\begin{lem}\label{kb2}
Suppose that $G_S$ has two S-cycle faces $F_1,F_2$ which lie locally
on the side of $T$ corresponding to the string $I_{1,2}$, such that
the circles $\partial F_1,\partial F_2$ are not isotopic in the
closed surface $T\cup I_{1,2}$. Then $T$ is generated by a once
punctured Klein bottle $P$ with $\partial P\subset I_{1,2}$, which
is essential whenever $(M,T_0)$ is not cabled and $M(\partial T)$ is
irreducible.
\end{lem}

\bpf
As observed above, the presence of S-cycle faces in $G_S$ implies
that $T$ is neutral, hence the surface $T\cup I_{1,2}$ is closed,
orientable, and of genus two. Since the circles $\partial
F_1,\partial F_2$ intersect the string $I_{1,2}$ each in one
spanning arc, and are disjoint and not isotopic in $T\cup I_{1,2}$,
compressing the surface $T\cup I_{1,2}$ in $M$ along the disks
$F_1,F_2$ produces a 2-sphere embedded in $M$, which bounds a 3-ball
in $M$ since $M$ is irreducible. It follows that $T$ separates $M$
into two components with closures $T^+,T^-$, so that if $T^+$ is the
component containing the string $I_{1,2}$ then $T^+$ is a genus two
handlebody with complete disk system $F_1,F_2$. Moreover, if $x,y$
are generators of $\pi_1(T^+)$ which are dual to $F_1,F_2$,
respectively, then, with some orientation convention, if $c$ is the
core of $I_{1,2}$ then $c$ represents the word $x^2y^2$ in
$\pi_1(T^+)$. As $c$ intersects each disk $F_1,F_2$ coherently in
two points, it is not hard to see that $c$ bounds a once punctured
Klein bottle $P$ in $T^+$ such that $T^+$ is homeomorphic to $N(P)$.

Finally, if $M(\partial T)$ is irreducible then $\wh{P}$ is
incompressible in $M(\partial T)$, so $P$ is incompressible in $M$
since $T_0$ is incompressible; and if $P$ boundary compresses in $M$
then it boundary compresses into a Moebius band, whence $(M,T_0)$ is
$(1,2)$-cabled. The lemma follows.
\epf

\section{Edge size}\label{neg}

In this section we will assume that $(T,\partial T)\subset (M,T_0)$
is a punctured torus with $t=|\partial T|\geq 1$ and $S$ a properly
embedded surface in $M$ which intersects $T$ in essential graphs
$G_S,G_T$, and establish bounds for the sizes of the edges in the
reduced graph $\bgs$, under suitable conditions. We denote the
vertices $S\cap T_0$ of $G_S$ by $u_i$'s, and the vertices of $G_T$
by $v_j$'s; notice that all edges in $G_S$ are internal.

\subsection{Positive edges}\label{possec}
A bound for the size of a positive edge of $\bgs$ can be easily
found.

\begin{lem}\label{pos}
Suppose $(M,T_0)$ is not cabled. If $t\geq 3$ and $\ove$ is a
positive edge of $\bgs$ then $|\ove|\leq t$, and if $|\ove|=t$ then
$t$ is even, the edge orbit of $\ove$ is a subgraph of $\bgt$
isomorphic to the graph of Fig.~\ref{n03} (thick edges only), and
some vertex of $\bgt$ has at most two incident positive nonloop
edges.
\end{lem}

\bpf
Let $t\geq 3$ and $\ove$ be a positive edge of $\bgs$ of size $\geq
t$, with consecutive edges $e_1,e_2,\dots,e_{t},e_{t+1},\dots$
labeled and running from $u_i$ to $u_{i'}$, as shown in
Fig.~\ref{n03}(a). The collection  $E=\{e_1,e_2,\dots,e_{t}\}$ then
induces a permutation $\sigma$ of the form $x\mapsto \alpha-x$, a
nontrivial involution, so the edge orbits of $E$ are a family of
disjoint cycles of length 2, which are nontrivial in $\wh{T}$ by
Lemma~\ref{mainref}(b), and hence the subgraph of $G_T$ generated by
these cycle edge orbits is isomorphic to the graph shown in
Fig.~\ref{n03}(b) (thick edges only). In particular, $t$ is even, so
$t\geq 4$, and there are $t/2$ such cycles. Consider now the the
edges $e_1,e_{\alpha-1}$, which form a cycle edge orbit of $E$ in
$G_T$ with vertices $v_1,v_{\alpha-1}$ of opposite parity. If
$|\ove|\geq t+1$ then, as the edge $e_{t+1}$ also has endpoints on
$v_1\cup v_{\alpha-1}$, it must lie in $T$ in the annular region
between the cycle formed by $e_1,e_{\alpha-1}$ and some other cycle
of $E$, which implies that $e_{t+1}$ is parallel to $e_1$ or
$e_{\alpha-1}$ in $T$, contradicting Lemma~\ref{mainref}(b) (see
Fig.~\ref{n03}(b)). Therefore $|\ove|\leq t$.

If $|\ove|=t$ then every vertex $v$ of $G_T$ belongs to a unique
cycle edge orbit $c(v)$ of $\ove$. Suppose that the vertices in
$c(v)$ are $v$ and $v'$. Then it is not hard to see from
Fig.~\ref{n03}(b) that $v$ can have at most two incident positive
nonloop edges of $\bgt$ on each side of the cycle $c(v)$; so if $v$
has at least three incident positive nonloop edges of $\bgt$, then
$v'$ can have at most one incident positive nonloop edge in $\bgt$
(see Fig.~\ref{n03}(b)).
\begin{figure}
\psfrag{(a)}{$(a)$}
\psfrag{(b)}{$(b)$}
\psfrag{e1}{$e_1$}
\psfrag{e2}{$e_2$}
\psfrag{et}{$e_t$}
\psfrag{et1}{$e_{t+1}$}
\psfrag{ek}{$e_{\alpha-1}$}
\psfrag{p}{$+$}
\psfrag{m}{$-$}
\psfrag{v}{$v$}
\psfrag{v'}{$v'$}
\psfrag{ui}{$u_i$}
\psfrag{ui'}{$u_{i'}$}
\Figw{5in}{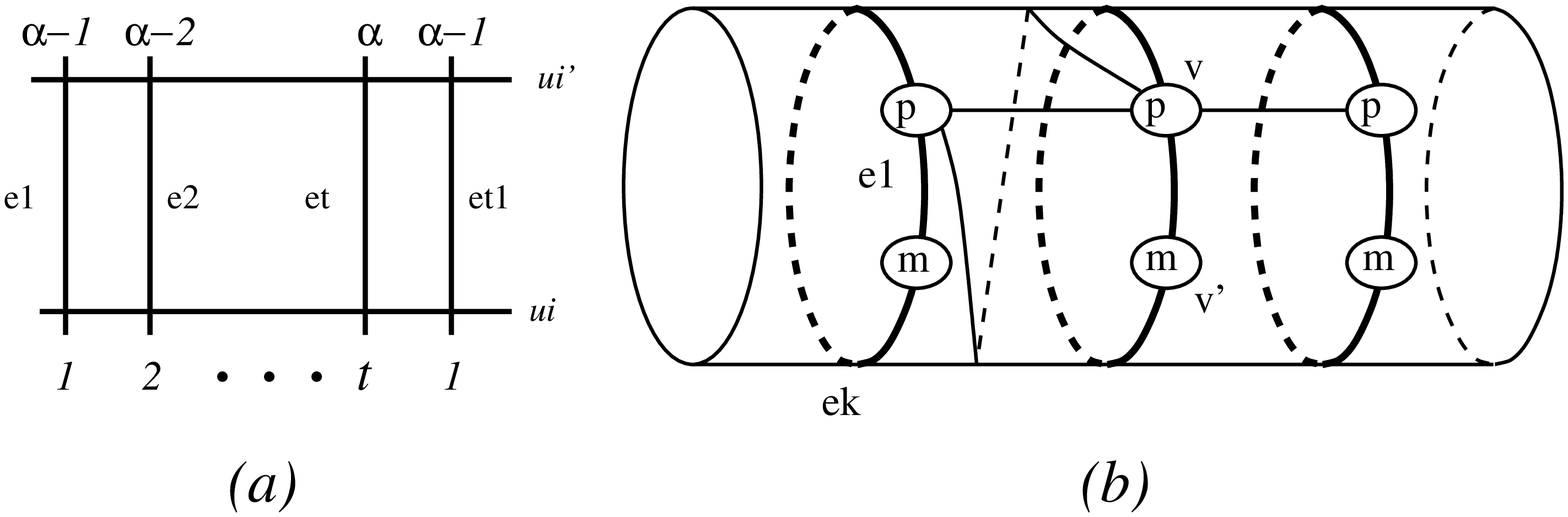}{}{n03}
\end{figure}
\epf

\subsection{Negative edges I}\label{negsec1}

The following fact is the starting point for our analysis of the
size of the negative edges in $\bgs$; its proof follows from
\cite[Lemma 2.8(2)]{tera6}, and we include it for the convenience of the
reader.

\begin{lem}\label{neg0}
Suppose $(M,T_0)$ is not cabled. If $t\geq 1$ and $\ove$ is a
negative edge of $\bgs$ with $|\ove|\geq t+1$, then $T$ is polarized
and any subcollection of $t$ consecutive edges in $\ove$ has exactly
one edge orbit. In particular, all disk faces of $G_S$ are even
sided.
\end{lem}

\bpf
Suppose $t\geq 1$ and there is a negative edge $\ove$ in $\bgs$ of
size $|\ove|\geq t+1$, with one endpoint in $u_i$ and the other in
$u_{i'}$. We may assume $e_1,\dots,e_t,e_{t+1},\dots$ are all the
edges in $\ove$, as shown in Fig.~\ref{n01}, oriented from $u_i$ to
$u_{i'}$.
\begin{figure}
\psfrag{ui}{$u_i$}
\psfrag{ui'}{$u_{i'}$}
\psfrag{e1}{$e_1$}
\psfrag{e2}{$e_2$}
\psfrag{e3}{$e_3$}
\psfrag{et}{$e_t$}
\psfrag{et1}{$e_{t+1}$}
\psfrag{et2}{$e_{t+2}$}
\psfrag{f1}{$F_1$}
\psfrag{f2}{$F_2$}
\psfrag{ft}{$F_t$}
\psfrag{ft1}{$F_{t+1}$}
\Figw{4.5in}{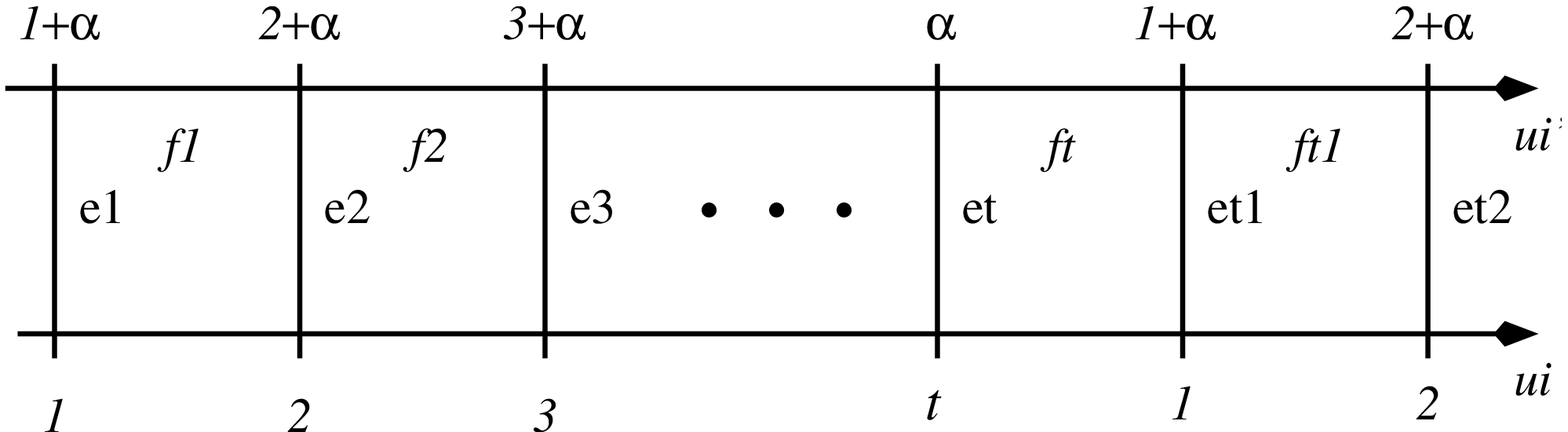}{}{n01}
\end{figure}

The collections of edges $E=\{e_1\dots,e_t\}$ and
$E'=\{e_2\dots,e_{t+1}\}$ induce the same permutation $\sigma$, of
the form $\sigma(x)=x+\alpha$ for some $0\leq\alpha<t$ (cf
\S\,\ref{orbits}), and
$\sigma$ has $n=\gcd(t,\alpha)$ orbits. By Lemma~\ref{mainref}(b),
in $G_T$, the edge orbits of each collection $E,E'$ are nontrivial
disjoint cycles and the edges $e_1$ and $e_{t+1}$ are not parallel.
Let $\gamma,\gamma'$ be the edge orbits of $E,E'$, that contain the
edges $e_1,e_{t+1}$, respectively. If $n\geq 2$ then the edge
$e_{t+1}$ is necessarily located in between two distinct edge orbits
of $E$, with both endpoints on the same side of the cycle $\gamma$
in $G_T$, as shown in Fig.~\ref{n02}. As the edges of $\gamma'$
coincide with those of $\gamma$, except for the edge $e_1$ which
gets replaced by $e_{t+1}$, it follows that the cycle $\gamma'$
bounds a disk in $\wh{T}$, contradicting Lemma~\ref{mainref}(b).
Therefore $n=1$, so $\sigma$, and hence $E$, have a single orbit,
and so $T$ is polarized; thus, by the parity rule, all edges in
$G_S$ are negative, from which it follows that any boundary
component of any face of $G_S$ has an even number of sides.
\begin{figure}
\psfrag{e1}{$e_1$}
\psfrag{et1}{$e_{t+1}$}
\psfrag{g}{$\gamma$}
\Figw{3in}{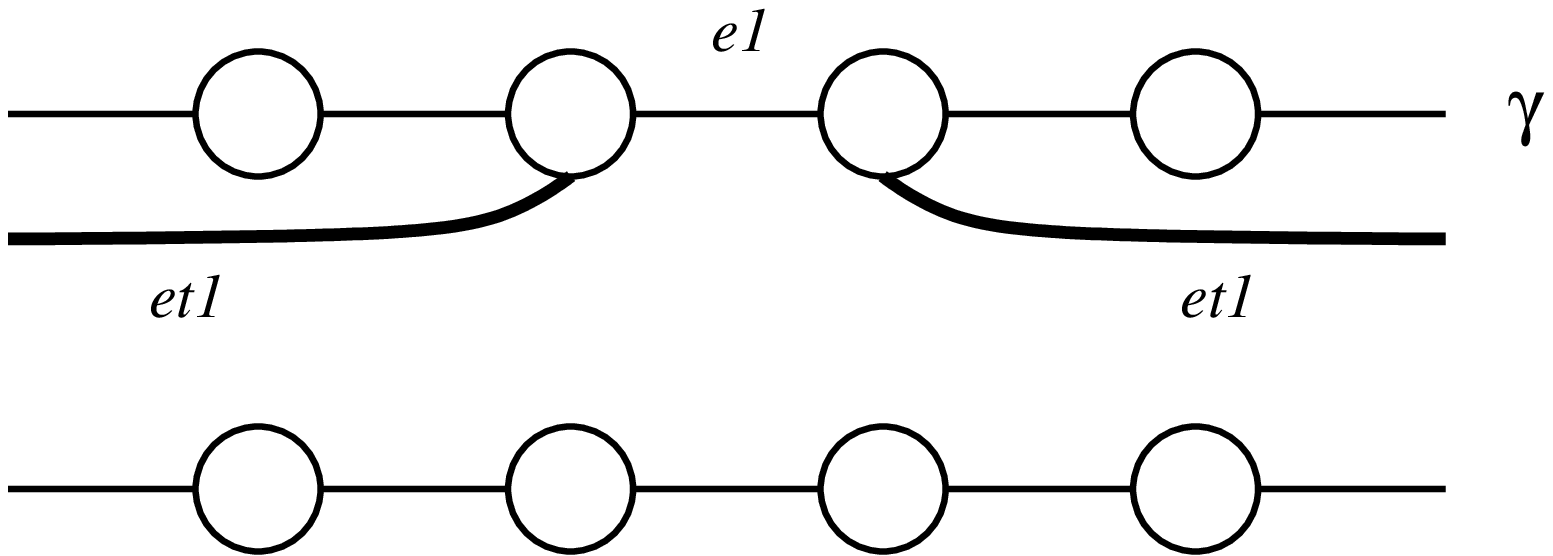}{}{n02}
\end{figure}
\epf

\subsection{Negative edges II: a construction}\label{negsec2}

For any properly embedded surface $\mc{F}$ in  a 3-manifold
$\mc{M}$, we will denote by $\mc{M}_{\mc{F}}=\mc{M}\setminus\intr
N(\mc{F})$ the manifold obtained by cutting $\mc{M}$ along $\mc{F}$;
if $\mc{F}$ is orientable then $N(\mc{F})=\mc{F}\times I$, where
$I=[0,1]$, and $\partial \mc{M}_{\mc{F}}$ contains two copies
$\mc{F}^0=\mc{F}\times 0,
\mc{F}^1=\mc{F}\times 1$ of $\mc{F}$.

Given any collection $E$ of mutually parallel, consecutive negative
edges of $G_S$ with $|E|\geq 2$, we define $M_{T,E}\subset M$ as the
submanifold obtained by cutting $M$ along the union of $T$ and the
bigon faces cobounded by the edges of $E$ in $G_S$. In this section
we will take a closer look at the manifolds $M_{T,E}$ constructed
with large enough collections $E$. Observe $M_T$ and $M_{T,E}$ are
irreducible manifolds.

Let $t\geq 1$ and $E=\{e_1,e_2,\dots,e_{t+1}\}$ be any collection of
$t+1$ mutually parallel, consecutive, negative edges in $G_S$,
running and oriented from the vertex $u_i$ to the vertex $u_{i'}$ of
$G_S$, and labeled as in Fig.~\ref{n01}.  By Lemma~\ref{neg1}, $T$
is polarized, hence nonseparating in $M$, so the permutation induced
by $E$ is of the form $x\mapsto x+\alpha\mod t$ with
$\gcd(t,\alpha)=1$. In what follows, for clarity, our figures will
sometimes be sketched to represent scenarios for large $t$, but the
arguments and constructions can be seen to hold for all $t\geq 1$.

It follows from the proof of Lemma~\ref{neg0} that the union of the
edges in $E$ form a subgraph of $G_T$ isomorphic to the graph shown
in Fig.~\ref{n04}, where $e_1,\dots,e_t$ are represented by the
horizontal edges and $e_{t+1}$ by the thicker edge. Moreover, if
$\gamma_1,\gamma_2$ are the cycle edge orbits of the collections
$\{e_1,\dots,e_t\}$ and $\{e_2,\dots,e_{t+1}\}$, respectively, then
$\gamma_1$ is the oriented cycle comprised of all the  horizontal
edges in Fig.~\ref{n04}, while $\gamma_2$ is obtained from
$\gamma_1$ by exchanging the edge $e_1$ with the edge $e_{t+1}$.
Hence $\Delta(\gamma_1,\gamma_2)=1$ in $\wh{T}$. The situation gets
somewhat simplified in the case $t=1$ from what is shown in
Fig.~\ref{n04}, which deals with the cases $t\geq 2$.
\begin{figure}
\psfrag{g}{$\gamma_1$}
\psfrag{e1}{$e_1$}
\psfrag{et1}{$e_{t+1}$}
\psfrag{mt}{$\mu_{t-\alpha}$}
\psfrag{m1}{$\mu_{1}$}
\psfrag{m1a}{$\mu_{1+\alpha}$}
\psfrag{N}{$\vec{N}$}
\Figw{3.5in}{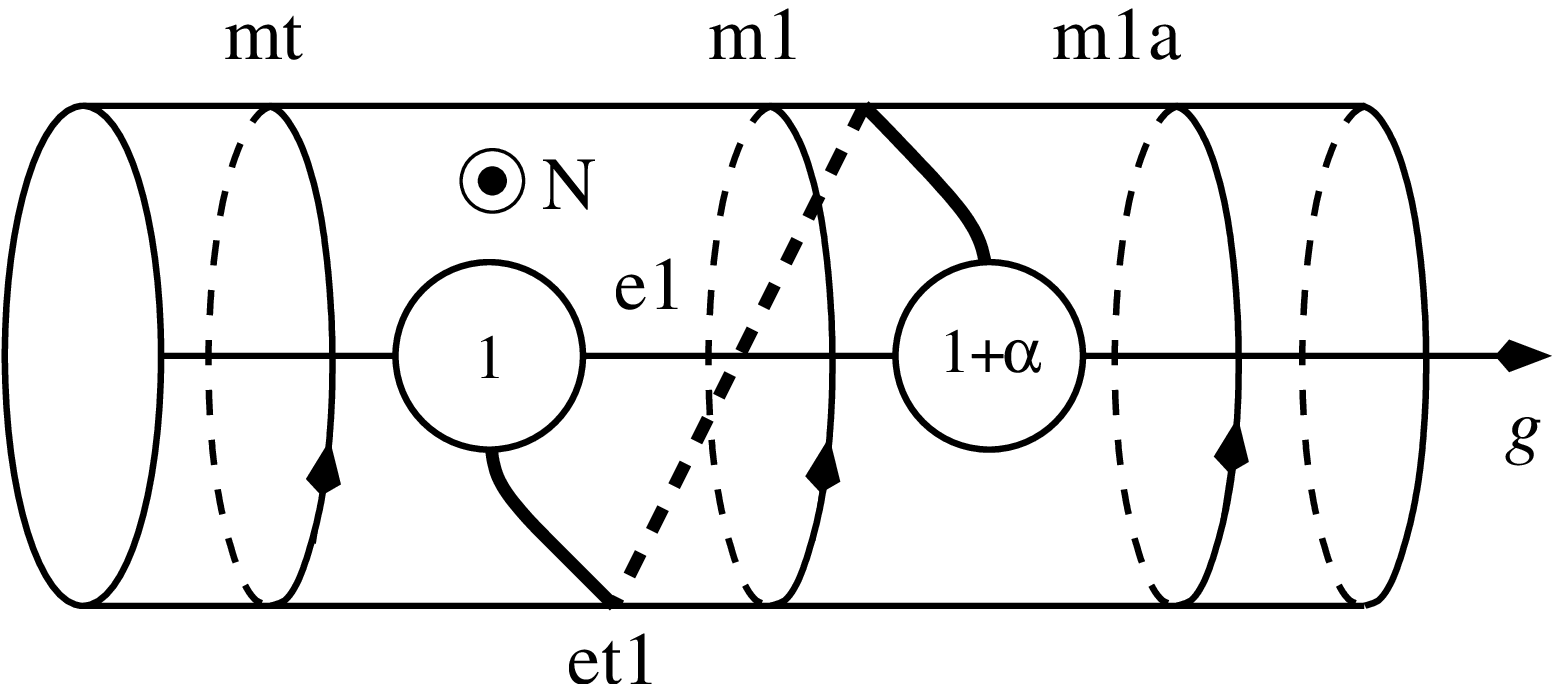}{}{n04}
\end{figure}

Fig.~\ref{n04} also shows a collection of $t$ oriented circles
$\mu_i$, $1\leq i\leq t$, each having the same slope in $\wh{T}$ as
the cycle $e_1\cup e_{t+1}$ and labeled by the vertex along
$\gamma_1$ that precedes (following the orientation of $\gamma_1$).
Notice that the cycles $\gamma_1$ and $\gamma_2$ can be obtained
from each other via one full Dehn twist on $T$ along $\mu_1$.

Each vertex $v_k$, edge $e_k$, and circle $\mu_k$ in $T$ splits into
two copies $v^1_k, e^1_k,\mu^1_k\subset T^1$, $v^2_k,
e^2_k,\mu^2_k\subset T^2$, with $v^1_k,v^2_k$, $e^1_k,e^2_k$, and
$\mu^1_k,\mu^2_k$ parallel in $N(T)=T\times I$ to $v_k$, $e_k$,
$\mu_k$, respectively; Fig.~\ref{n11} shows such parallelism for
$e_k,e^1_k,e^2_k$.

Let $\psi:T^1\rightarrow T^2$ be the gluing homeomorphism that
produces $M$ out of $M_T$. We will orient $e^1_k,e^2_k$ and
$\mu^1_k,\mu^2_k$ in the same direction as $e_k,\mu_k$,
respectively, via the parallelism $N(T)=T\times I$, so that
$\psi(e^1_k)=e^2_k$ and $\psi(\mu^1_k)=\mu^2_k$, preserving
orientations.

The edges $e^1_1,e^1_2,\dots,e^1_t$ form a cycle in $T^1$ which is
parallel in $N(T)$ to the cycle $\gamma_1$, while
$e^2_2,e^2_3,\dots,e^2_{t+1}$ form a cycle in $T^2$ parallel in
$N(T)$ to $\gamma_2$. We will denote these cycles by
$\gamma_1^1\subset T^1$ and $\gamma_2^2\subset T^2$, respectively.
Thus, while the cycles $\gamma_1$ and $\gamma_2$ intersect in $T$,
the cutting process along $T$ `separates' them into the disjoint
cycles $\gamma_1^1,\gamma_2^2$, with $\psi(\gamma^1_1)=\gamma_2^2$.
\begin{figure}
\psfrag{ui}{$u_i$}
\psfrag{ui'}{$u_{i'}$}
\psfrag{e1}{$e_1$}
\psfrag{e2}{$e_2$}
\psfrag{e3}{$e_3$}
\psfrag{et}{$e_t$}
\psfrag{et1}{$e_{t+1}$}
\psfrag{t1}{$T^1$}
\psfrag{t2}{$T^2$}
\psfrag{f1}{$F'_1$}
\psfrag{f2}{$F'_2$}
\psfrag{ft}{$F'_t$}
\psfrag{e21}{$e^2_1$}
\psfrag{e11}{$e^1_1$}
\psfrag{e2t1}{$e^2_{t+1}$}
\psfrag{e1t1}{$e^1_{t+1}$}
\psfrag{I'12}{$I'_{1,2}$}
\psfrag{I't1}{$I'_{t,1}$}
\Figw{4.5in}{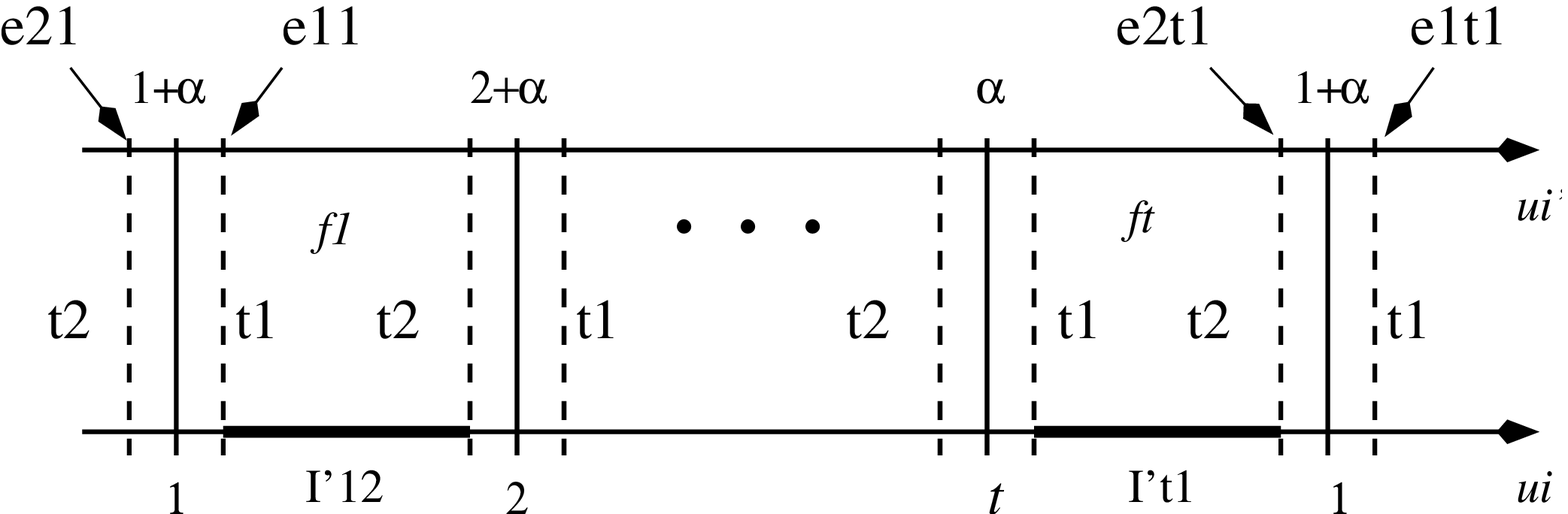}{}{n11}
\end{figure}

For each string $I_{k,k+1}$ of $T$, we will call the annulus
$I'_{k,k+1}=I_{k,k+1}\cap M_T\subset\partial M_T$ a {\it string of
$M_T$}. Observe that the union of $T^1,T^2$, and the strings of
$M_T$ is one of the boundary components of $M_T$, of genus $t+1$.

Consider now the bigons $F_1,F_2,\dots,F_t$ of $G_S$ cobounded by
the edges of $E$, as shown in Fig.~\ref{n01}. We call the disks
$F'_k=F_k\cap M_T\subset M_T$, $1\leq k\leq t$, the {\it faces} of
$E$ in $M_T$; these faces have corners along the strings
$I'_{k,k+1}$ and are properly embedded in $M_T$. For $1\leq k\leq
t$, $\partial F'_k$ consists of four segments: one corner along the
string $I'_{k,k+1}$, one corner along $I'_{k+\alpha,k+\alpha+1}$,
and the two edges $e^1_k\subset T^1$ and $e^2_{k+1}\subset T^2$ (see
Fig.~\ref{n11}). Since $|E|=t+1$, along each vertex $u_i,u_{i'}$
each string of $T$ appears exactly once among the corners of the
bigon faces $F_k$. Thus, each string $I'_{k,k+1}$ of $M_T$ has
exactly two corners coming from all the faces $F'_j$ in $E$, and
these two corners cut $I'_{k,k+1}$ into two rectangular pieces,
which we denote by $J_{k,k+1}, L_{k,k+1}$.

It follows that the faces $F'_k$ are embedded in $M_T$ as shown in
Fig.~\ref{n05}(a). To determine the location of the edges
$e^2_1\subset T^2$ and $e^1_{t+1}\subset T^1$, consider the normal
vector $\vec{N}$ to $T$ indicated in Fig.~\ref{n04} by the tip of an
arrow $\odot$, and orient $T^1,T^2$ via normal vectors
$\vec{N}^1,\vec{N}^2$, respectively, such that
$\vec{N}^1=\vec{N}^2=\vec{N}$ after identifying $T^1$ with $T^2$;
these vectors are indicated in Fig.~\ref{n05}(a), and we will use
them to identify the right hand and left hand sides of the cycles
$\gamma\subset T,\gamma^1_1\subset T^1,\gamma^2_2\subset T^2$
consistently. Since the oriented edge $e_{t+1}$ has initial and
terminal endpoints on the right and left hand sides of the oriented
cycle $\gamma$, respectively, the endpoints of the edge
$e^i_{t+1}\subset T^i$ must behave the same way relative to the
oriented cycle of edges $e^i_1\cup e^i_2\cup\dots\cup e^i_t\subset
T^i$ for $i=1,2$. Therefore the edges $e^2_1\subset T^2,
\ e^1_{t+1}\subset T^1$ must be embedded as shown in Fig.~\ref{n05}(a) (up
to Dehn twists in the annuli $T^i\setminus\gamma_i^i$), and hence
$\mu^1_1,\mu^2_1$ must then be embedded in $T^1,T^2$ as shown in
Fig.~\ref{n05}(a).
\begin{figure}
\psfrag{g1}{$\gamma_1^1$}
\psfrag{g2}{$\gamma_2^2$}
\psfrag{f'1}{$F'_1$}
\psfrag{f't}{$F'_t$}
\psfrag{I12}{$I'_{1,2}$}
\psfrag{Ia12}{$I'_{1+\alpha,2+\alpha}$}
\psfrag{It1}{$I'_{t,1}$}
\psfrag{Ia1}{$I'_{\alpha,\alpha+1}$}
\psfrag{B}{$B$}
\psfrag{m11}{$\mu^1_1$}
\psfrag{m21}{$\mu^2_1$}
\psfrag{N1}{$\vec{N}^1$}
\psfrag{N2}{$\vec{N}^2$}
\psfrag{e21}{$e^2_1$}
\psfrag{e22}{$e^2_2$}
\psfrag{e11}{$e^1_1$}
\psfrag{e2t1}{$e^2_{t+1}$}
\psfrag{e1t1}{$e^1_{t+1}$}
\psfrag{e1t}{$e^1_{t}$}
\psfrag{T1}{$T^1$}
\psfrag{T2}{$T^2$}
\psfrag{J12}{$J'_{1,2}$}
\psfrag{J1a2a}{$J'_{1+\alpha,2+\alpha}$}
\psfrag{L12}{$L'_{1,2}$}
\psfrag{L1a2a}{$L'_{1+\alpha,2+\alpha}$}
\psfrag{L}{$L'_{k,k+1}$}
\psfrag{AT}{$A_E$}
\psfrag{A'T}{$A'_E$}
\psfrag{(a)}{$(a)$}
\psfrag{(b)}{$(b)$}
\Figw{5in}{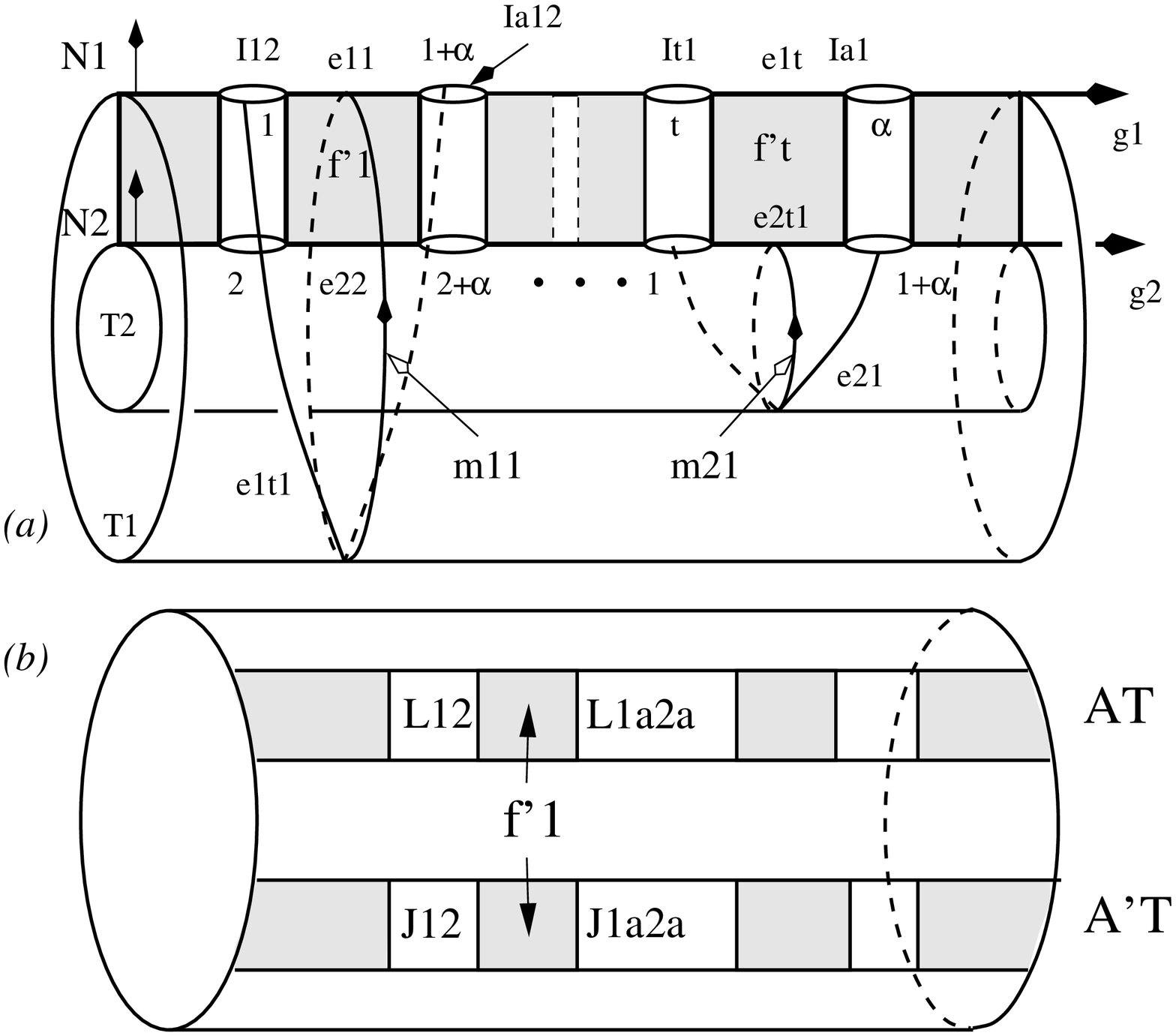}{}{n05}
\end{figure}

Cutting $M_T$ along the faces $F'_k$ produces the irreducible
submanifold $M_{T,E}\subset M$, which has a distinguished torus
boundary component $R_E$ that contains all the rectangles
$J_{k,k+1}, L_{k,k+1}$ and two copies of each face $F'_k$. The union
of all these pieces forms two disjoint nontrivial annuli
$A_E,A'_E\subset R_E$; relabeling if necessary, we may assume that
$A_E$ contains all the rectangles $J_{k,k+1}$, while $A'_E$ contains
the $L_{k,k+1}$'s (see Fig.~\ref{n05}(b)).

So, if $\wh{M}_T$ is the manifold obtained by cutting $M(\partial
T)$ along $\wh{T}$, it is not hard to see that $\wh{M}_T$ can be
obtained from $M_{T,E}$ by identifying $A_E$ with $A'_E$ in such a
way that all pairs of rectangles $J_{k,k+1}$ and $L_{k,k+1}$ match.

A first approximation to the structures of $M,M_T$, and $M(\partial
T)$ is given in our next result.

\begin{lem}\label{neg1}
Suppose that $t\geq 1$ and $\ove$ is a negative edge of $\bgs$ of
size $|\ove|\geq t+1$.
\ben
\ita If $T$ is incompressible in $M$, $|\ove|=t+1$, and the torus
$R_{\ove}\subset\partial M_{T,\ove}$ compresses in $M_{T,\ove}$,
then $\partial M=T_0$ and $\wh{M}_T$ is a Seifert fibered space over
the annulus with at most one singular fiber;

\itb
if $|\ove|\geq t+2$ then $M_T\approx T\times I$, so $\partial
M=T_0$, $T$ is incompressible in $M$, and $M(\partial T)$ is an
irreducible torus bundle over the circle with fiber $\wh{T}$.
\een
\end{lem}

\bpf
For part (a), let $c$ be the core of the annulus $A_{\ove}\subset
R_{\ove}$. If $D$ is a compression disk for $R_{\ove}$ in
$M_{T,\ove}$ then, as $T$ is incompressible and $M_{T,\ove}$ is
irreducible, we must have $d=\Delta(\partial D,c)\geq 1$ and
$M_{T,\ove}$ a solid torus. Hence $\partial M=T_0$ and $\wh{M}_T$ is
a Seifert fibered space over the annulus with at most one singular
fiber, of index $d$.

For part (b), let $\{e_1,e_2,\dots,e_{t+1},e_{t+2}\}$ be a
collection of $t+2$ consecutive edges in $\ove$, with edges and
bigons labeled as in Fig.~\ref{n01}, and consider the manifold
$M_{T,E}$ corresponding to the family of edges
$E=\{e_1,e_2,\dots,e_{t+1}\}$. As the edge $e^1_{t+1}$ is not
parallel in $T^1$ to any of the edges of the cycle
$\gamma_1^1\subset T^1$ (see Fig.~\ref{n05}(a)), the disk face
$F'_{t+1}=F_{t+1}\cap M_T$ is not parallel in $M_T$ to any of the
disks $F'_k=F_k\cap M_T$ for $1\leq k\leq t$, and hence $F'_{t+1}$
is necessarily embedded in $M_T$ as shown in Fig.~\ref{n07} (with
$\partial F'_{t+1}$ the union of the thicker edges
$e^1_{t+1},e^2_{t+2}$ and corners). It follows that $F'_{t+1}$,
which also lies in $M_{T,E}$, intersects each annulus $A_E,A'_E$
transversely in one spanning arc. Therefore, by the argument of part
(a), $\wh{M}_T$ is a Seifert fibered space over the annulus with no
singular fibers, so $\wh{M}_T\approx\wh{T}\times I$, from which it
follows that $M_T\approx T\times I$ and $M(\partial T)$ is an
irreducible torus bundle over the circle with incompressible fiber
$\wh{T}$.
\begin{figure}
\psfrag{g1}{$\gamma_1^1$}
\psfrag{g2}{$\gamma_2^2$}
\psfrag{f'1}{$F'_1$}
\psfrag{e22}{$e^2_2$}
\psfrag{e11}{$e^1_1$}
\psfrag{e2t2}{$e^2_{t+2}$}
\psfrag{e1t1}{$e^1_{t+1}$}
\psfrag{T1}{$T^1$}
\psfrag{T2}{$T^2$}
\Figw{3.5in}{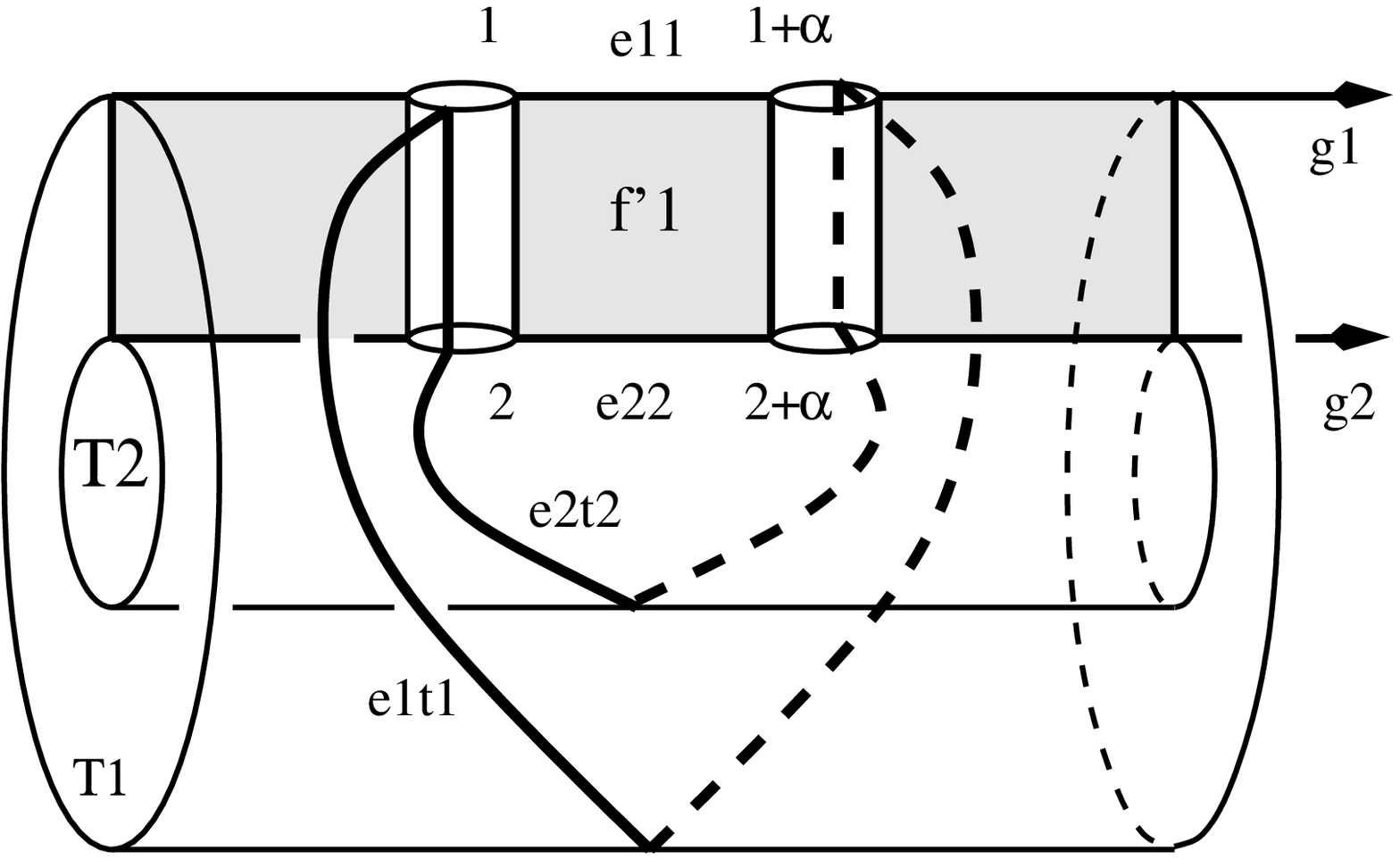}{}{n07}
\end{figure}
\epf

\subsection{Negative edges III}\label{negsec3}

In this section we will further assume that $G_S$ has at least $t+2$
mutually parallel negative edges, no two of which are parallel in
$T$, and determine the structure of $M$ under these conditions.

We will use the following definitions. Let $P$ be an oriented pair
of pants with boundary components $\mu_0,\mu_1,\mu_2$, each given
the induced orientation from $P$, and consider the manifold $P\times
S^1$, where $P$ is identified with some fixed copy $P\times \{*\}$
in $P\times S^1$. We orient the manifold $P\times S^1$ via a product
orientation, so that the circles $\{*\}\times S^1$ all follow the
direction of an orientation normal vector of $P$ in $P\times S^1$;
the boundary tori components $T_i=\mu_i\times S^1, \ i=0,1,2$, are
then oriented by an outside pointing normal vector $\vec{N}_i$.

Let $\phi:T_1\rightarrow T_2$ be an orientation reversing
homeomorphism such that
\begin{equation}\label{eq0}
\phi(\mu_1)=-\mu_2,
\end{equation}
where $-\mu_2$ is the circle $\mu_2$ with the opposite orientation.
Then the quotient manifold $P\times S^1/\phi$ is orientable,
irreducible, and has incompressible boundary the torus $T_0$. Also,
under the quotient map $P\times S^1\rightarrow P\times S^1/\phi$,
$P$ gives rise to a once punctured torus $T_P$ in $P\times S^1/\phi$
with boundary slope $\mu_0$, and the tori $T_1,T_2$ give rise a
closed, nonseparating, incompressible torus $T'_P\subset P\times
S^1/\phi$ which intersects $T_P$ transversely in a single circle
corresponding to $\mu_1=-\mu_2$.

Consider the arcs $h_i\subset P, \ 0\leq i\leq 5$, shown in
Fig.~\ref{n22}; these arcs give rise to essential annuli $h_i\times
S^1\subset P\times S^1$, which are the unique (up to isotopy)
properly embedded essential annuli in $P\times S^1$; in fact, the
annulus $h_0\times S^1$ is the unique essential surface in $P\times
S^1$ with boundary on $T_0$ (a similar statement holds for
$h_4\times S^1$ and $h_5\times S^1$). In particular, the pair
$(P\times S^1,T_0)$ is not cabled.
\begin{figure}
\psfrag{l0}{$\lambda_0$}
\psfrag{l1}{$\lambda_1$}
\psfrag{l2}{$\lambda_2$}
\psfrag{m1}{$\mu_1$}
\psfrag{m0}{$\mu_0$}
\psfrag{m2}{$\mu_2$}
\psfrag{h1}{$h_1$}
\psfrag{h0}{$h_0$}
\psfrag{h2}{$h_2$}
\psfrag{h3}{$h_3$}
\psfrag{h4}{$h_4$}
\psfrag{h5}{$h_5$}
\psfrag{T0}{$T_0$}
\psfrag{T1}{$T_1$}
\psfrag{T2}{$T_2$}
\Figw{4.5in}{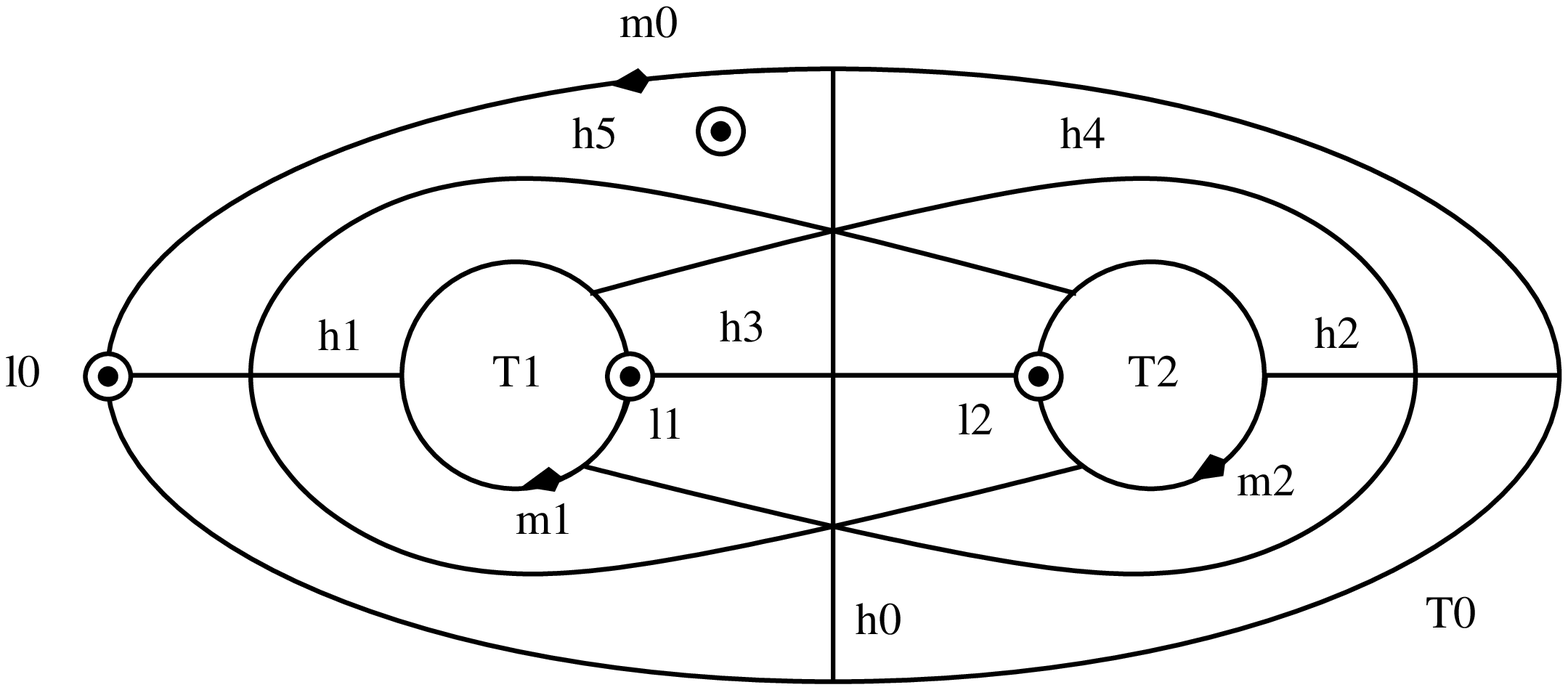}{}{n22}
\end{figure}
The boundary components of these annuli correspond to three slopes
$\lambda_i\subset T_i, \ i=0,1,2$; these are the unique slopes on
the $T_i$'s arising from any Seifert fibration on $P\times S^1$. We
will orient all the circles $\mu_i,\lambda_i$, $i=0,1,2$, as shown
in Fig.~\ref{n22}, where the tips of arrows $\odot$ indicate
directions of tangent/normal vectors, and each circle $\mu_i$ is
also labeled by the torus $T_i$ that contains it; thus, the
$\lambda_i$'s have the same orientation as the fibers $\{x\}\times
S^1$. The first homology group $H_1(P\times S^1)$ (with integer
coefficients) is then freely generated by $\mu_1,\mu_2,\lambda_0$,
and the following relations hold:
\begin{equation}\label{eq1}
\lambda_0=\lambda_1=\lambda_2, \ \mu_0+\mu_1+\mu_2=0.
\end{equation}

The above orientation scheme allows us to define intersection
numbers between two oriented  circles $c,c'$ in any boundary torus
$T_i$ of $P\times S^1$, by requiring that $c\cdot c'$ be positive at
a point $x\in c\cap c'$ of transverse intersection iff the tangent
vectors $\vec{v},\vec{v}\,'$ to $c,c'$ at $x$, respectively, yield
an orientation triple $(\vec{v},\vec{v}\,',\vec{N}_i)$ of $P\times
S^1$ at $x$. With this convention, the fact that $\phi$ is
orientation reversing can be restated as follows:
\begin{equation}\label{eq2}
\text{for any two oriented circles }c,c'\subset T_1, \
\phi(c)\cdot \phi(c')=-c\cdot c'.
\end{equation}

Notice that $\mu_1\cdot\lambda_1=+1$ and $\mu_2\cdot\lambda_2=+1$.
Relative to these orientation frames $\mu_1,\lambda_1$ of $H_1(T_1)$
and $\mu_2,\lambda_2$ of $H_1(T_2)$, we can write
$\phi(\lambda_1)=m\mu_2+r\lambda_2$ in $H_1(T_2)$ for some
relatively prime integers $m,r$, and then it follows from
(\ref{eq0}), (\ref{eq2}), and  $\mu_1\cdot\lambda_1=+1=
\mu_2\cdot\lambda_2$ that $r=+1$, so
\begin{equation}\label{eq3}
\phi(\lambda_1)=m\mu_2+\lambda_2\text{ \ in \ }H_1(T_2).
\end{equation}

The homeomorphism $\phi$ is determined up to isotopy by its action
on first homology, which, relative to the orientation frames
$\mu_1,\lambda_1$ and $\mu_2,\lambda_2$, is given by the matrix
$[\phi]=\left(\begin{smallmatrix}
-1 & m \\
0 & 1
\end{smallmatrix}\right)$;
so we will also denote the manifold $P\times S^1/\phi$ by $P\times
S^1/[m]$, and use the notation $$(P\times
S^1/[m],T_P,T'_P,T_0,\mu_0,\lambda_0)$$ to stress the presence of
the specific objects $T_P,T'_P,T_0,\mu_0,\lambda_0\subset P\times
S^1/\phi$ constructed above. Since the quotient manifolds $P\times
S^1/\phi$ and $P\times S^1/\phi^{-1}$ are homeomorphic, and
$[\phi^{-1}]=\left(\begin{smallmatrix}
-1 & -m \\
0 & 1
\end{smallmatrix}\right)$,
switching the roles of $T_1,T_2$ in $P\times S^1/\phi^{-1}$ gives
rise to a homeomorphism $P\times S^1/[m]\approx P\times S^1/[-m]$,
and so we may assume that $m\geq 0$. Finally, we identify the cut
manifold $(P\times S^1/[m])_{T'_P}$ with $P\times S^1$.

The main result of this section can now be stated as follows:

\begin{prop}\label{prop1}
Let $(T,\partial T)\subset (M,T_0)$ be a $\mc{K}$-incompressible
torus and $S\subset M$ a surface which intersects $T$ in essential
graphs $G_S,G_T$. Set $t=|\partial T|\geq 1$, and suppose that $G_S$
has at least $t+2$ mutually parallel, consecutive negative edges, no
two of which are parallel in $T$. If $t=1$ then $M$ is the exterior
of the trefoil knot, while if $t\geq 2$ then $M=(P\times
S^1/[m],T_P,T'_P,T_0,\mu_0,\lambda_0)$ with $T_P$ having the same
boundary slope $\mu_0$ as $T$, and the following hold:
\ben
\ita
$(M,T_0)$ is not cabled;

\itb
$M(\alpha)$ is irreducible and toroidal for any slope
$\alpha\neq\lambda_0$, and $M(\lambda_0)\approx S^1\times S^2\# L$
for some closed 3-manifold $L$ of genus $\leq 1$;

\itc
$M=P\times S^1/[m]$ contains a punctured $\mc{K}$-incompressible
torus with boundary slope $\alpha\neq\mu_0$ iff $m=1,2,4$ and
$\alpha$ is the slope of $\mu_0-(4/m)\lambda_0$; in such case,
$\Delta(\alpha,\mu_0)=4/m=1,2,4$ and $M$ also contains an essential
$q$-punctured Klein bottle of boundary slope $\alpha$, where
$(m,q)=(1,1), \ (2,1)$, or $(4,2)$.\qed
\een
\end{prop}

Proposition~\ref{prop1} follows immediately from Lemmas~\ref{t1},
\ref{neg2}, \ref{neg3}, and \ref{neg4} below. In what follows we
will also use the notation of \S\ref{negsec2}; as usual, we may draw
some figures as if $t$ were large only for clarity.

\begin{lem}\label{t1}
If $t=1$ then $M$ is the exterior of the trefoil knot.
\end{lem}

\bpf
Let $e_1,e_2,e_3$ be three distinct mutually parallel, consecutive
edges in $G_S$ which are not parallel in $G_T$, running from $u_i$
to $u_{i'}$, and let $F_1,F_2$ be the bigon disk faces they cobound
in $G_S$, as shown in Fig.~\ref{n21}(a). The graph $\bgt$ is then
isomorphic to the graph shown in Fig.~\ref{n18}(a).
\begin{figure}
\psfrag{e11}{$e^1_1$}
\psfrag{e22}{$e^2_2$}
\psfrag{e21}{$e^2_1$}
\psfrag{e13}{$e^1_3$}
\psfrag{e12}{$e^1_2$}
\psfrag{e23}{$e^2_3$}
\psfrag{f'1}{$F'_1$}
\psfrag{pd}{$F'_2$}
\psfrag{T1}{$T^1$}
\psfrag{T2}{$T^2$}
\psfrag{e1}{$e_1$}
\psfrag{e2}{$e_2$}
\psfrag{e3}{$e_3$}
\psfrag{u1}{$v_1$}
\psfrag{vj}{$u_i$}
\psfrag{vj'}{$u_{i'}$}
\psfrag{(a)}{$(a)$}
\psfrag{(b)}{$(b)$}
\psfrag{(c)}{$(c)$}
\psfrag{1}{$1$}
\Figw{4.5in}{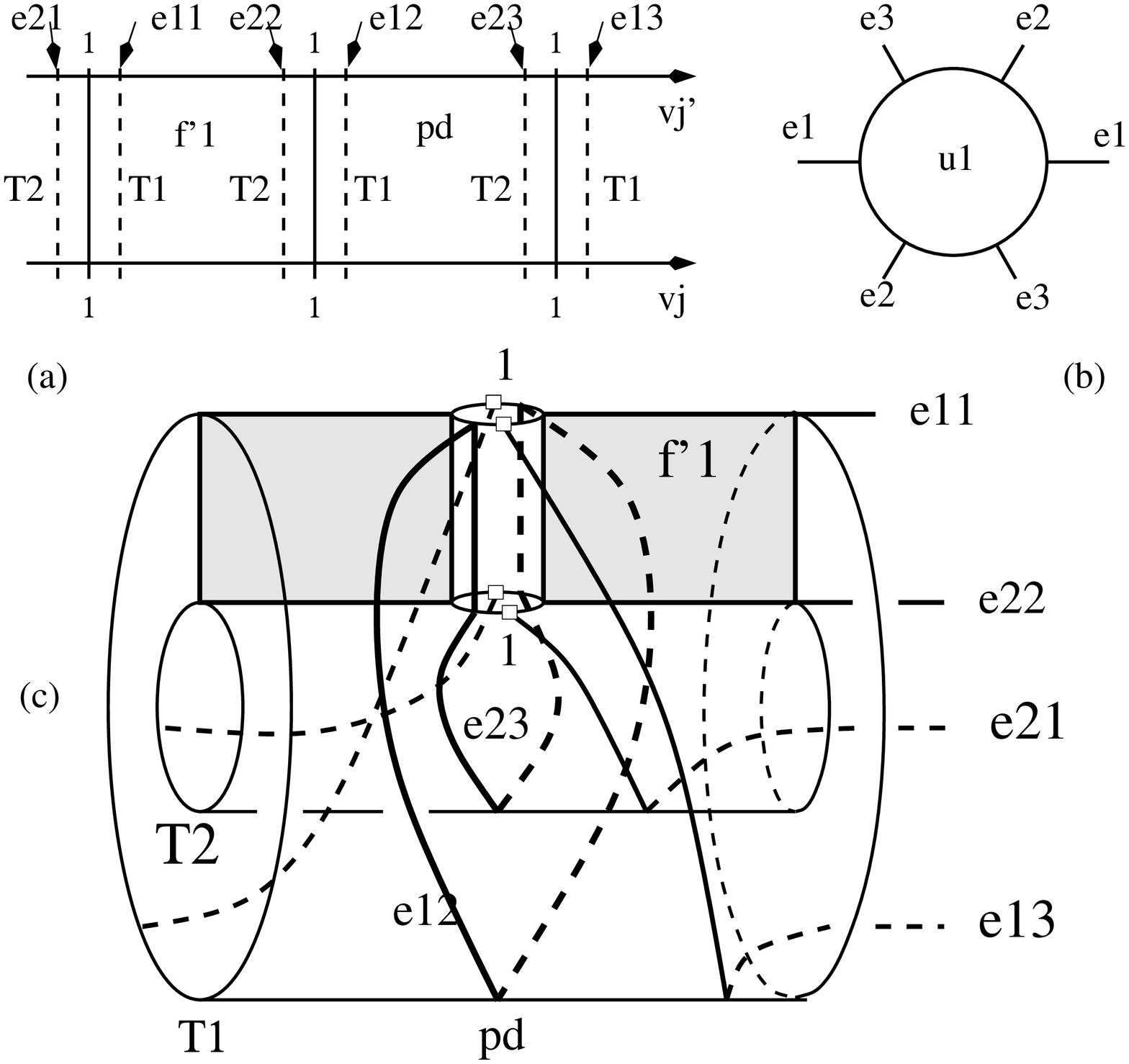}{}{n21}
\end{figure}

As in the proof of Lemma~\ref{neg1}, the edges $e^1_1, e^2_2$ and
$e^1_2,e^2_3$ must lie in $T^1,T^2$ as shown in Fig.~\ref{n21}(c),
cobounding the faces $F'_1,F'_2$ of $M_T$, respectively. To locate
the edges $e^2_1\subset T^2$ and $e^1_3\subset T^1$, first observe
that locally the edges $e_1,e_2,e_3$ produce the pattern around the
vertex $v_1$ of $G_T$ shown in Fig.~\ref{n21}(b), say with the edges
$e_1,e_2,e_3$ repeating cyclically twice around $v_1$, in that
order, so that exactly the same pattern must be present around each
copy $v^1_1\subset T^1, \ v^2_1\subset T^2$ of $v_1$. Therefore, if
$e^1_3$ is embedded in $T^1$ as shown in Fig.~\ref{n21}(c), then
$e^2_1$ must be embedded in $T^2$ as shown in Fig.~\ref{n21}(c), and
so these two edges $e^1_3,e^2_1$, along with two spanning arcs on
the string $I'_{1,1}$ of $M_T$, cobound a rectangular disk $D$ in
$M_T$ disjoint from $F'_1\cup F'_2$ (but not necessarily from $S$).

Let $B\subset M$ be the surface obtained from the union of
$F'_1,F'_2,D$, after identifying $T^1$ with $T^2$ in $M_T$ via
$\psi$ so that $e^1_k=e^2_k$ for $k=1,2,3$. Then $B$ is either an
annulus or a Moebius band in $(M,T_0)$ which intersects $T$ in
essential graphs consisting of exactly 3 edges, so that the graph
$G_{T,B}=B\cap T\subset T$ has two triangle faces $C_1,C_2$ (see
Fig.~\ref{n18}(a)).

If $B$ is an annulus then $\Delta(\partial B,\partial T)=3$, $B$ is
neutral by the parity rule, and the faces $C_1$ and $C_2$ locally
lie on opposite sides of $B$. Hence cutting the irreducible manifold
$M$ along $B$ produces two solid tori $V_1,V_2$ with corresponding
meridian disks $C_1,C_2$. Since all edges of the graph $G_{T,B}$ are
positive, each disk $C_1,C_2$ intersects the annulus $B$ coherently
and transversely in 3 spanning arcs, and so $M=V_1\cup_B V_2$ is
homeomorphic to a Seifert fibered space with base a disk and two
singular fibers of indices $3,3$. However, $T$ is the union of the
two meridian disks $C_1$ and $C_2$ along the 3 edges of $B\cap T$,
and it is not hard to see that any such union produces a pair of
pants in $M$, not a once punctured torus. Therefore $B$ must be a
Moebius band, with $\Delta(\partial B,\partial T)=6$, and cutting
$M$ along $B$ produces a solid torus with meridian disk either
triangle face $C_1$ or $C_2$ of $G_{T,B}$. Since all edges of the
graph $G_{T,B}$ are positive, $M$ must be homeomorphic to a Seifert
fibered space with base a disk and two singular fibers of indices
$2,3$, which is the trefoil knot exterior.
\epf

The rest of this section is devoted to the cases $t\geq 2$.

\begin{lem}\label{neg2}
If $t\geq 2$ then $M=(P\times S^1/[m],T_P,T'_P,\mu_0,\lambda_0)$ for
some integer $m$, with $T_P$ and $T$ having the same boundary slope
$\mu_0$.
\end{lem}

\bpf
Suppose that $t\geq 2$ and $\ove$ is a negative edge of $\bgs$ with
$|\ove|\geq t+2$. We assume that
$e_1,\dots,e_t,e_{t+1},e_{t+2},\dots$ are all the edges in $\ove$,
labeled as in Fig.~\ref{n01}, and oriented from $u_i$ to $u_{i'}$.

Let $\psi:T^1\rightarrow T^2$ be the gluing map that produces $M$
out of $M_T$. As in the proof of Lemma~\ref{neg1}, the face
$F'_{t+1}=F_{t+1}\cap M_T$ is properly embedded in $M_T$ with
boundary as shown in Fig.~\ref{n07}, and $M_T=T\times I$.

Consider now the oriented circles $\mu_1,\dots,\mu_t$ embedded in
$T$ as shown in Fig.~\ref{n04}. Recall each circle $\mu_k$ is
labeled by the vertex of $T$ that precedes it along the oriented
cycle $\gamma_1$ generated by $e_1,\dots e_t$ in $T$, and that
$\mu_k$ splits into two copies $\mu^1_k\subset T^1$ and
$\mu^2_k\subset T^2$, which are oriented in the same direction as
$\mu_k$ within $N(T)$; thus, all circles $\mu^1_k,\mu^2_k$ are
coherently oriented in $T^1,T^2$. From Figs.~\ref{n05}(a) and
\ref{n07} and the fact that $e^2_1,e^2_{t+2}$ are disjoint and
nonparallel in $T^2$, it follows that all circles $\mu^1_k\subset
T^1$ and $\mu^2_k\subset T^2$ are embedded as shown in
Fig.~\ref{n05-2}.

Therefore, the faces $F'_1$ and $F'_{t+1}$ can be isotoped in $M_T$
to construct an annulus $A_1\subset M_T$ with boundary the circles
$\mu^1_1\cup\mu^2_2$, which under their given orientations remain
coherently oriented relative to $A_1$. Via the product structure
$M_T=T\times I$, it is not hard to see that each pair of circles
$\mu^1_k,\mu^2_{k+1}$ cobounds such an annulus $A_k\subset M_T$ for
$1\leq k\leq t$, with the oriented circles $\mu^1_k,\mu^2_{k+1}$
coherently oriented relative to $A_k$; these annuli $A_k$ can be
taken to be mutually disjoint and $I$-fibered in $M_T=T\times I$.
Since $\psi(\mu^1_k)=\mu^2_k$ (preserving orientations), the union
$A_1\cup A_2\dots\cup A_t$ yields a closed nonseparating torus $T''$
in $M$, on which the circles $\mu_1,\mu_2,\dots,\mu_t$ appear
consecutively, coherently oriented, and in this order.
\begin{figure}
\psfrag{g1}{$\gamma_1^1$}
\psfrag{g2}{$\gamma_2^2$}
\psfrag{e22}{$e^2_2$}
\psfrag{e11}{$e^1_1$}
\psfrag{e21}{$e^2_1$}
\psfrag{u11}{$\mu^1_1$}
\psfrag{u22}{$\mu^2_2$}
\psfrag{u1k}{$\mu^1_k$}
\psfrag{u2k1}{$\mu^2_{k+1}$}
\psfrag{N1}{$\vec{N}^1$}
\psfrag{N2}{$\vec{N}^2$}
\psfrag{e2t2}{$e^2_{t+2}$}
\psfrag{e1t1}{$e^1_{t+1}$}
\psfrag{e1t}{$e^1_{k}$}
\psfrag{T1}{$T^1$}
\psfrag{T2}{$T^2$}
\Figw{5in}{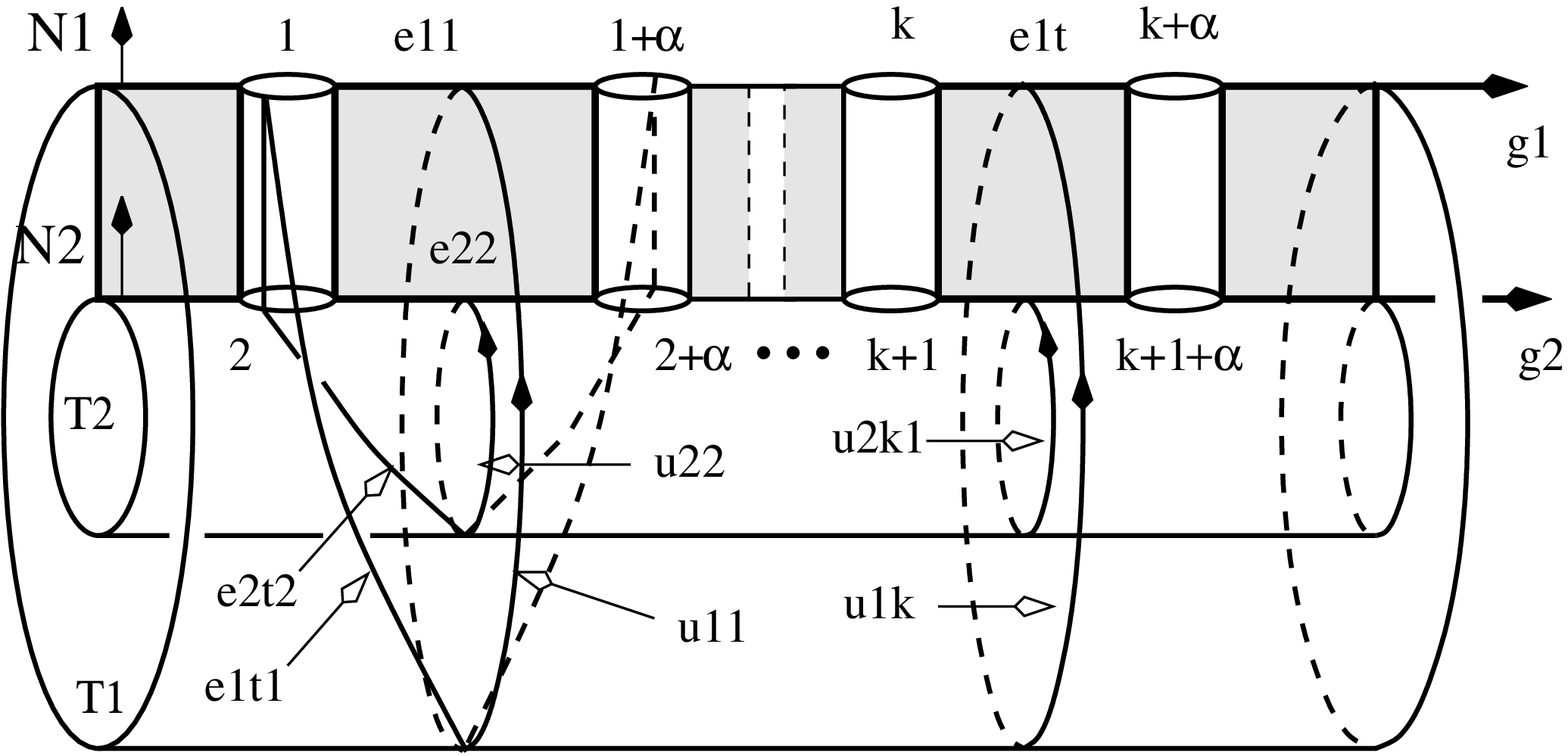}{}{n05-2}
\end{figure}

Thus, the region in $M_T$ cobounded by any pair $A_i,A_j$ of
consecutive annuli has a product structure of the form
$P_{i,j}\times I$, where $P_{i,j}$ is the pants region cobounded by
the boundary circles of $A_i,A_j$ in $T^1$; since $M_{T''}$ is the
union of these regions $P_{i,j}\times I$, glued along their pant
boundary pieces $P_{i,j}\times 0, \ P_{i,j}\times 1$ via the map
$\psi$, it follows that $M_{T''}$ has a product structure of the
form $P\times S^1$, where $P$ is any of the pants $P_{i,j}$.

Now, in $M$, $P$ has one boundary component $\partial_0 P$ on $T_0$,
of the same slope as $\partial T$, while the other two boundary
components $\partial_1 P,\partial_2 P$ lie on $T''$ and are
disjoint. From the point of view of $T''$, the circles $\partial_1
P,\partial_2 P$ coincide with two of the circles $\mu_k\subset T''$,
whose given orientations are coherent along $T''$; with such
orientations, the circles $\mu_k$ are then also coherently oriented
relative to $P$. Therefore $P$ can be isotoped in $M$ so that
$\partial_1P=\partial_2P$ on $T''$, giving rise to a once punctured
torus $T_P$ in $M$ of the same boundary slope as $T$ which
intersects $T''$ in a circle of the same slope as the $\mu_k$'s.

It follows that $M$ is a manifold of the form $(P\times
S^1/[m],T_P,T'_P,T_0,\mu_0,\lambda_0)$, with $T_P=T'$ and
$T'_P=T''$.
\epf

\begin{lem}\label{neg3}
If $M=(P\times S^1/[m],T_P,T'_P,T_0,\mu_0,\lambda_0)$ then $(M,T_0)$
is not cabled, $M(\lambda_0)\approx S^1\times S^2\# L$ for some
closed 3-manifold $L$ of genus $\leq 1$, and $M(\alpha)$ is
irreducible with $T'_P$ incompressible in $M(\alpha)$ for
$\alpha\neq\lambda_0$.
\end{lem}

\bpf
Write $M=P\times S^1/\phi$ with $\phi$ the gluing map
$\phi:T_1\rightarrow T_2$. Clearly, for any slope
$\alpha\neq\lambda_0$, $(P\times S^1)(\alpha)$ is an irreducible
Seifert fibered space over an annulus with at most one singular
fiber. Therefore the tori $T_1\cup T_2=\partial (P\times
S^1)(\alpha)$ are incompressible in $(P\times S^1)(\alpha)$, so
$M(\alpha)=(P\times S^1)(\alpha)/\phi$ is irreducible and hence the
nonseparating torus $T'_P$ is incompressible in $M(\alpha)$.

Consider now the manifold $M(\lambda_0)$. Let $A$ be the
nonseparating and neutral annulus $h_0\times S^1\subset P\times
S^1$, of boundary slope $\lambda_0$. Then $\wh{A}$ is a
nonseparating 2-sphere in $M(\lambda_0)$ disjoint from the
nonseparating torus $T'_P$. Observe that $T'_P$ compresses in
$M(\lambda_0)$, on both sides, via the disks generated by the annuli
$h_1\times S^1$ and $h_2\times S^1$ of $P\times S^1$, whose
boundaries are the circles $\lambda_1\subset T_1$ and
$\lambda_2\subset T_2$. Thus, cutting $M(\lambda_0)$ along
$\wh{A}\cup T'_P$ yields two once punctured solid tori $V_1,V_2$
with torus boundary components $T_1,T_2$ and meridian disks of
boundary slopes $\lambda_1\subset T_1,\lambda_2\subset T_2$,
respectively. Gluing $V_1$ to $V_2$ along $T_1,T_2$ via $\phi$ then
produces a twice punctured manifold $L^-$; since
$|\phi(\lambda_1)\cdot\lambda_2|=m$, identifying the two spherical
boundary components of $L^-$ via $\phi$ produces the manifold
$M(\lambda_0)= S^1\times S^2\# L$, where $L=S^1\times S^2, \ S^3$,
or a lens space for $m=0,1$, or $m\geq 2$, respectively.

Finally, suppose that $(M,T_0)$ is cabled, with essential cabling
annulus $A'$. Then $A'$ is separating and hence neutral; as the
annulus $A$ is also neutral, by the parity rule $A$ and $A'$ must
have the same boundary slope $\lambda_0\subset T_0$. Since $A$ is
the unique essential surface in $P\times S^1$ with boundary on
$T_0$, it follows that, after isotoping $A'$ in $M$ so as to
intersect $T'_P$ transversely and minimally, we must have $|A'\cap
T'_P|>0$. Thus $A'\cap P\times S^1$ is a collection of essential
annular components, each of which must then be isotopic to some
annulus $h_i\times S^1\subset P\times S^1$, $1\leq i\leq 5$. It is
not hard to see that two such components must be isotopic to $h
_1\times S^1$ and $h_2\times S^1$, which implies that
$\phi(\lambda_1)=\lambda_2$, whence $m=0$ and so
$M(\lambda_0)=S^1\times S^2\# S^1\times S^2$ by the argument above,
contradicting the fact that $A'$ being a cabling annulus implies
that $M(\partial A')=M(\lambda_0)$ has a lens space connected
summand.
\epf

\begin{lem}\label{neg4}
$M=P\times S^1/[m]$ contains a punctured $\mc{K}$-incompressible
torus with boundary slope $\alpha\neq\mu_0$ iff $m=1,2,4$ and
$\alpha$ is the slope of $\mu_0-(4/m)\lambda_0$; in such case,
$\Delta(\alpha,\mu_0)=4/m=1,2,4$ and $M$ also contains an essential
$q$-punctured Klein bottle of boundary slope $\alpha$, where
$(m,q)=(1,1), \ (2,1)$, or $(4,2)$.
\end{lem}

\bpf
Let $R$ be a punctured essential torus or Klein bottle in $M$. Since
the only connected essential surface in $P\times S^1$ with boundary
on $T_0$ is the annulus $h_0\times S^1$, after isotoping $R$ in $M$
so that it intersects $T'_P$ transversely and minimally, we must
have $|R\cap T'_P|> 0$, with each circle component of $R\cap T'_P$
nontrivial in both $R$ and $T'_P$ and each component of $R'=R\cap
P\times S^1$ essential in $P\times S^1$. Isotoping $R'$, we may
assume that $R'$ and $P$ intersect transversely in essential graphs.

\setcounter{claim*}{0}
\begin{claim*}\label{c00}
$\alpha\neq\lambda_0$, so $M(\alpha)$ is irreducible and each
component of $R\cap T'_P$ is nontrivial in both $\wh{R}$ and $T'_P$.
\end{claim*}

Clearly, there is some edge $x$ in the essential graph $R'\cap
P\subset P$ for which at least one endpoint lies on $T_0$, ie, $x$
is isotopic in $P$ to $h_0,h_1$, or $h_2$. If $\alpha=\lambda_0$
then, as $R'$ is essential, the annulus $A'=x\times S^1$ can be
isotoped in $P\times S^1$ so as to be disjoint from $R'$, which
implies that $R'$ lies in the cut manifold $(P\times S^1)_{A'}$,
where it is incompressible. As $(P\times S^1)_{A'}$ consists of one
or two copies of $\text{(torus)}\times I$'s, it follows that $R'$
must be a union of annuli, and hence that $R$ is an annulus, which
is not the case. Therefore $\alpha\neq\lambda_0$ and hence, by
Lemma~\ref{neg3}, $M(\alpha)$ is irreducible and the nonseparating
torus $T'_P$ is incompressible in $M(\alpha)$, whence each component
of $R\cap T'_P$ must be nontrivial in both $\wh{R}$ and $T'_P$.
\qed(Claim~\ref{c00})

Thus, by Lemma~\ref{neg3} and Claim~\ref{c00}, $(M,T_0)$ is not
cabled and $M(\alpha)$ is irreducible. Let $Q$ be any component of
$R\cap P\times S^1$; by Claim~\ref{c00}, $Q$ is an essential
punctured annulus with two boundary components
$\partial_1Q,\partial_2 Q$ in $T_1\cup T_2$, and without loss of
generality we may assume that $q=|\partial Q\cap T_0|>0$. If
$\partial_1Q\cup\partial_2 Q\subset T_1$ or
$\partial_1Q\cup\partial_2 Q\subset T_2$ then $Q$ boundary
compresses in $P\times S^1$ relative to $T_0$ via the annulus
$h_2\times S^1$ or $h_1\times S^1$, respectively, contradicting the
fact that $R$ is essential in $M$; thus we may assume that
$\partial_1Q\subset T_1$ and $\partial Q_2\subset T_2$.

\begin{claim*}\label{c4}
$\Delta(\alpha,\lambda_0)=1$ and all components of $\partial Q\cap
T_0$ are coherently oriented in $T_0$.
\end{claim*}

Isotope $Q$ so that it intersects the annuli $(h_0\cup h_1\cup
h_2)\times S^1\subset P\times S^1$ transversely in essential graphs;
then for $i=0,1,2$ each graph $h_i\times S^1\cap Q\subset h_i\times
S^1$ consists of $\Delta(\alpha,\lambda_0)\cdot q\geq q$ parallel
edges, all of which are internal and negative for $i=0$. The union
of any $q$ consecutive edges in the graphs $h_i\times S^1\cap
Q\subset h_i\times S^1$ for $i=1,2$ produce a subgraph in $Q$ of the
type shown in Fig.~\ref{n23} (vertical thin edges). Therefore, any
$q$ consecutive edges of the graph  $h_0\times S^1\cap Q\subset
h_0\times S^1$ necessarily lie in $Q$ like the thick horizontal
edges shown in Fig.~\ref{n23}, so any edge of the graph $h_0\times
S^1\cap Q\subset h_0\times S^1$ is parallel in $Q$ to some
horizontal edge of Fig.~\ref{n23}. Since the pair $(P\times
S^1,T_0)$ is not cabled, it follows from Lemma~\ref{mainref}(b) that
we must have $\Delta(\alpha,\lambda_0)\cdot q=q$, so
$\Delta(\alpha,\lambda_0)=1$, in which case the edges of the graph
$h_0\times S^1\cap Q\subset h_0\times S^1$, all of which are
negative, form a single cycle in $Q$, which implies that all the
components of $\partial Q\cap T_0$ are coherently oriented in $T_0$.
\qed(Claim~\ref{c4})
\begin{figure}
\psfrag{b1q}{$\partial_1Q$}
\psfrag{b2q}{$\partial_2Q$}
\psfrag{(a)}{$(a)$}
\psfrag{(b)}{$(b)$}
\psfrag{m1}{$\mu_1$}
\psfrag{m0}{$\mu_0$}
\psfrag{m2}{$\mu_2$}
\psfrag{T1}{$T_1$}
\psfrag{T2}{$T_2$}
\Figw{2.7in}{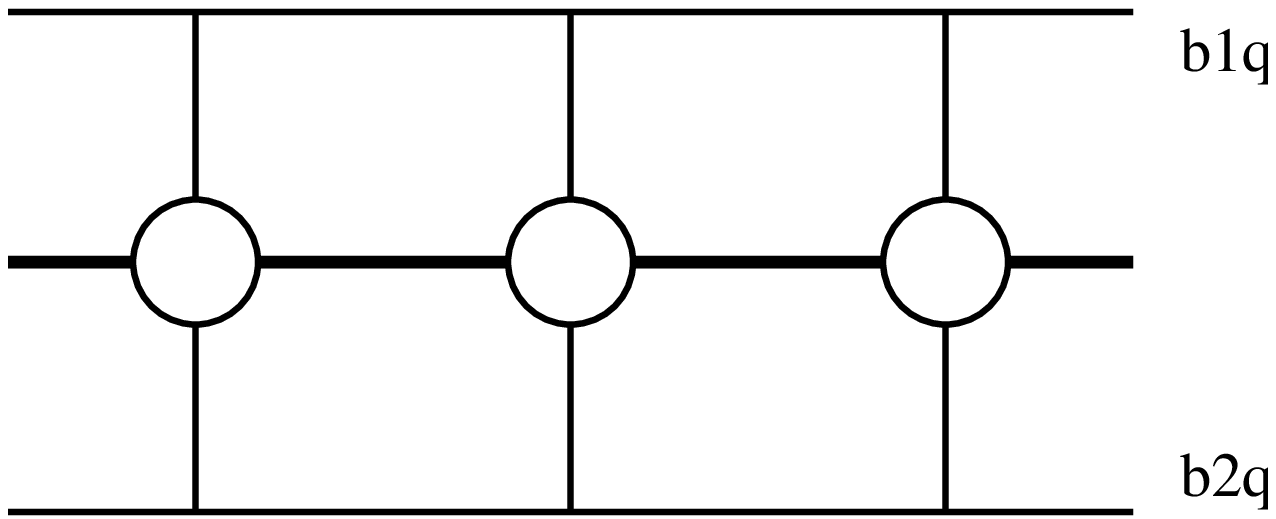}{}{n23}
\end{figure}

We now select the orientation on $Q$ which induces the orientation
on any component $c$ of $\partial Q\cap T_0$ such that $c\cdot
\lambda_0>0$. Since $\Delta(\alpha,\lambda_0)=1$ by Claim~\ref{c4}, we can
write $\alpha=\mu_0+b_0\lambda_0$, where
$|b_0|=\Delta(\alpha,\mu_0)$, and so
\begin{align*}
\partial Q\cap T_0& =q\alpha=q(\mu_0+b_0\lambda_0)\text{ in } H_1(T_0),\\
\partial_1 Q & =a_1\mu_1+b_1\lambda_1\text{ in } H_1(T_1), \ \
\text{ and } \ \ \\
\partial_2 Q & =a_2\mu_2+b_2\lambda_2\text{ in } H_1(T_2),
\end{align*}
for some pairs $a_1,b_1$ and $a_2,b_2$ of relatively prime integers.
Since $\partial Q=0$ in $H_1(P\times S^1)$, by (\ref{eq1}) we have
that
\begin{align*}
0=\partial Q\cap T_0+\partial_1Q+\partial_2Q =
(a_1-q)\mu_1+(a_2-q)\mu_2+(qb_0+b_1+b_2)\lambda_0,
\end{align*}
and hence $a_1=a_2=q$ and $qb_0=-(b_1+b_2)$, so that
\begin{equation}\label{eq4}
\partial_1 Q  =q\mu_1+b_1\lambda_1\text{ in } H_1(T_1)\text{ \ \ and \ \ }
\partial_2 Q  =q\mu_2+b_2\lambda_2\text{ in } H_1(T_2).
\end{equation}

Observe that $\phi(\partial_1Q)=\pm\partial_2 Q$ in $H_1(T_2)$ since
$\phi$ maps the circle $\partial_1Q\subset T_1$ onto a circle in
$T_2$ of the same slope as $\partial Q_2$.

\begin{claim*}\label{c5}
$\phi(\partial_1Q)=+\partial_2Q$ in $H_1(T_2)$, $(m,q)=(1,1), \
(2,1)$, or $(4,2)$, and $\alpha=\mu_0-(4/m)\lambda_0$.
\end{claim*}

We have $\phi(\partial_1Q)=\ve\partial_2 Q$ in $H_1(T_2)$ for some
$\ve\in\{\pm 1\}$; from (\ref{eq0}), (\ref{eq3}), and (\ref{eq4}),
it follows that
$$
-q\mu_2+b_1(m\mu_2+\lambda_2)=\ve(q\mu_2+b_2\lambda_2)
\text{ in }H_1(T_2),
$$
and hence that $b_1\, m=(1+\ve)q$ and $b_1=\ve b_2$. If $\ve=-1$
then $b_1=-b_2$ and so $qb_0=-(b_1+b_2)=0$; but then $b_0=0$, whence
$\alpha=\mu_0$, which is not the case. Hence $\ve=+1$, so
$b_1=b_2=2q/m$ and $b_0=-4/m$, and so
$\Delta(\alpha,\mu_0)=|b_0|=4/m\geq 1$. In particular, $m=1,2,4$ and
$\alpha=\mu_0-(4/m)\lambda_0$, and since $a_1=q$ and $b_1=2q/m$ are
relatively prime integers we must have $(m,q)=(1,1), \ (2,1)$, or
$(4,2)$.\qed(Claim~\ref{c5})

Therefore, for each pair $(m,q)$ listed in Claim~\ref{c5}, $Q$ is
$q$-punctured annulus which can be isotoped in $P\times S^1$ so that
$\phi(\partial_1 Q)=\partial_2Q$ in $T_2$, giving rise to a
$q$-punctured Klein bottle $Q'$ in $M=P\times S^1/[m]$ with boundary
slope $\alpha=\mu_0-(4/m)\lambda_0$.

Now, since $M(\alpha)$ is irreducible, the closed Klein bottle
$\wh{Q}'$ is necessarily incompressible in $M(\alpha)$. So, if $Q'$
is not essential in $M$ then a compression or boundary compression
of $Q'$ gives rise to either a $(q-1)$-punctured Moebius band $B$ in
$(M,T_0)$ or a closed Klein bottle $R''$ in $M$. In the first case,
$\wh{B}$ is a projective plane in the irreducible manifold
$M(\alpha)$, which implies that $M(\alpha)$ is homeomorphic to
$RP^3$, contradicting the fact that $M(\alpha)$ is a toroidal
manifold for $\alpha\neq\lambda_0$. And in the second case, the
closed Klein bottle $R''$ must be incompressible in the irreducible
manifold $M$, whence $R''$ can be isotoped to intersect $T'_P$
transversely and minimally, so that $|R''\cap T'_P|>0$ and $R''\cap
P\times S^1$ consists of annuli, all of which are essential in
$P\times S^1$; since any such annulus must then be isotopic to one
of the annuli $h_i\times S^1$, $i=3,4,5$, it follows that $m=0$,
which is not the case. Therefore $Q'$ is essential in $M$.

Conversely, let $(m,q)$ be one of the pairs $(1,1), \ (2,1), \
(4,2)$, and let $\alpha=\mu_0-(4/m)\lambda_0$; then a punctured
annulus $Q$ can be constructed in $P\times S^1$ with $q$ punctures
in $T_0$ of slope $\alpha$ and one puncture in $T_i$ of slope
$q\mu_i+(2q/m)\lambda_i$ for $i=1,2$, by homologically summing, in a
suitable way, $q$ copies of $P$ and $2q/m$ copies of each annulus
$h_1\times S^1, \ h_2\times S^1$. Since any homeomorphism
$\phi:T_1\rightarrow T_2$ that homologically maps $\mu_1$ onto
$-\mu_2$ and $\lambda_1$ onto $m\mu_2+\lambda_2$ also maps
$q\mu_1+(2q/m)\lambda_1$ onto $q\mu_2+(2q/m)\lambda_2$, the lemma
follows.
\epf

\section{Boundary slopes of $\mc{K}$-incompressible tori}
\label{delta6}

In this section we assume that $(M,T_0)$ is not cabled and that
$(F_1,\partial F_1)$ and $(F_2,\partial F_2)$ are
$\mc{K}$-incompressible tori in $(M,T_0)$ with boundary slopes at
distance $\Delta\geq 6$ and essential graphs of intersection; by
Lemma~\ref{mainref}(c), both Dehn filled manifolds $M(r_1)$ and
$M(r_2)$ are irreducible,

We will use the generic notation $\{S,T\}=\{F_1,F_2\}$, $s=|\partial
S|$, and $t=|\partial T|$, and denote the vertices of $G_S$ by
$u_i$'s and those of $G_T$ by $v_j$'s.

By Proposition~\ref{prop1}, for $t\geq 1$, any negative edge in
$G_S$ has size at most $t+1$. This last bound can be improved a bit
in many cases, given that $\Delta\geq 6$, as shown below.

\begin{lem}\label{lesst}
If $\Delta\geq 6$ and $t\geq 3$ then $\Delta=6$ and, in $\bgs$,
$\deg\equiv 6$ and any edge has size $t$.
\end{lem}

\bpf
By Lemma~\ref{pos}, any positive edge of $\bgs$ has size at most
$t$.

Suppose there is a negative edge $\ove$ in $\bgs$ of size $t+1$. By
Lemma~\ref{neg0}, any disk face of $\bgs$ is even sided, and so
$\bgs$ has a vertex $u_i$ of degree at most 4 by Lemma~\ref{3v}(b).
If $u_i$ has $p$ positive and $n$ negative local edges in $\bgs$
then $\dbgs(u_i)=p+n\leq 4$ and so the degree of $u_i$ in $G_S$
satisfies the relations
$$6t\leq \Delta\cdot t=\dgs(u_i)\leq p\cdot t+n\cdot (t+1)
=(p+n)t+n\leq 4t+4,
$$
whence $t\leq 2$, which is not the case.

Therefore any edge of $\bgs$ has size at most $t$, so if $u$ is any
vertex of $\bgs$ with $p'$ positive and $n'$ negative local edges,
then again the degree of $u$ in $G_S$ satisfies the relations
$$6t\leq \Delta\cdot t=\dgs(u)\leq p'\cdot t+n'\cdot t=(p'+n')t;$$
thus $\dbgs(u)=p'+n'\geq 6$, and hence $\deg\equiv 6$ in $\bgs$ by
Lemma~\ref{3v}(a). Since equality must then hold throughout the
above relations, it follows that $\Delta=6$ and each edge of $\bgs$
has size $t$.
\epf

The {\it jumping number} of the graphs $G_S$ and $G_T$ was
introduced in \cite[\S 2]{gordon5}. For $\Delta=6$ the jumping
number is one, which means that if the $\Delta$ points of
intersection between the circles $\partial_i S$ ($=u_i$) and
$\partial_j T$ ($=v_j$) are labeled consecutively as
$x_1,x_2,\dots,x_{\Delta}$ around $\partial_i S$, then these points
appear consecutively around $\partial_j T$ in the same order
$x_1,x_2,\dots,x_{\Delta}$ when read in some direction. We will
refer to this corresponding distribution of labels around the
vertices of $G_S$ and $G_T$ as the {\it jumping number one
condition,} or {\it JN1 condition} for short.

\subsection{The generic cases $s,t\geq 3$}

By Lemma~\ref{lesst}, $\Delta=6$ and, in $\bgs,\bgt$, $\deg\equiv 6$
and all edges have size $t,s$, respectively; in particular, for any
label $1\leq j\leq t$ ($1\leq i\leq s$), each vertex $w$ of $G_S$
($G_T$, resp.) has 6 local edges with label $j$ ($i$, resp.) at $w$,
which give rise to the 6 local edges around $w$ in $\bgs$ ($\bgt$,
resp.). The JN1 condition now implies that if
$e_1,e_2,e_3,e_4,e_5,e_6$ are the local edges with label $j$ at
$u_i$, as shown in Fig.~\ref{n14b}, then these local edges appear
with label $i$ around $v_j$ as shown in Fig.~\ref{n14b}, up to
reflection about a diameter of $v_j$; and, by the parity rule, any
local edge around $v_j$ has the opposite sign of the corresponding
local edge around $u_i$.
\begin{figure}
\psfrag{e1}{$e_1$}
\psfrag{e2}{$e_2$}
\psfrag{e3}{$e_3$}
\psfrag{e4}{$e_4$}
\psfrag{e5}{$e_5$}
\psfrag{e6}{$e_6$}
\psfrag{ui}{$u_i$}
\psfrag{vj}{$v_j$}
\psfrag{i}{$i$}
\psfrag{j}{$j$}
\Figw{3in}{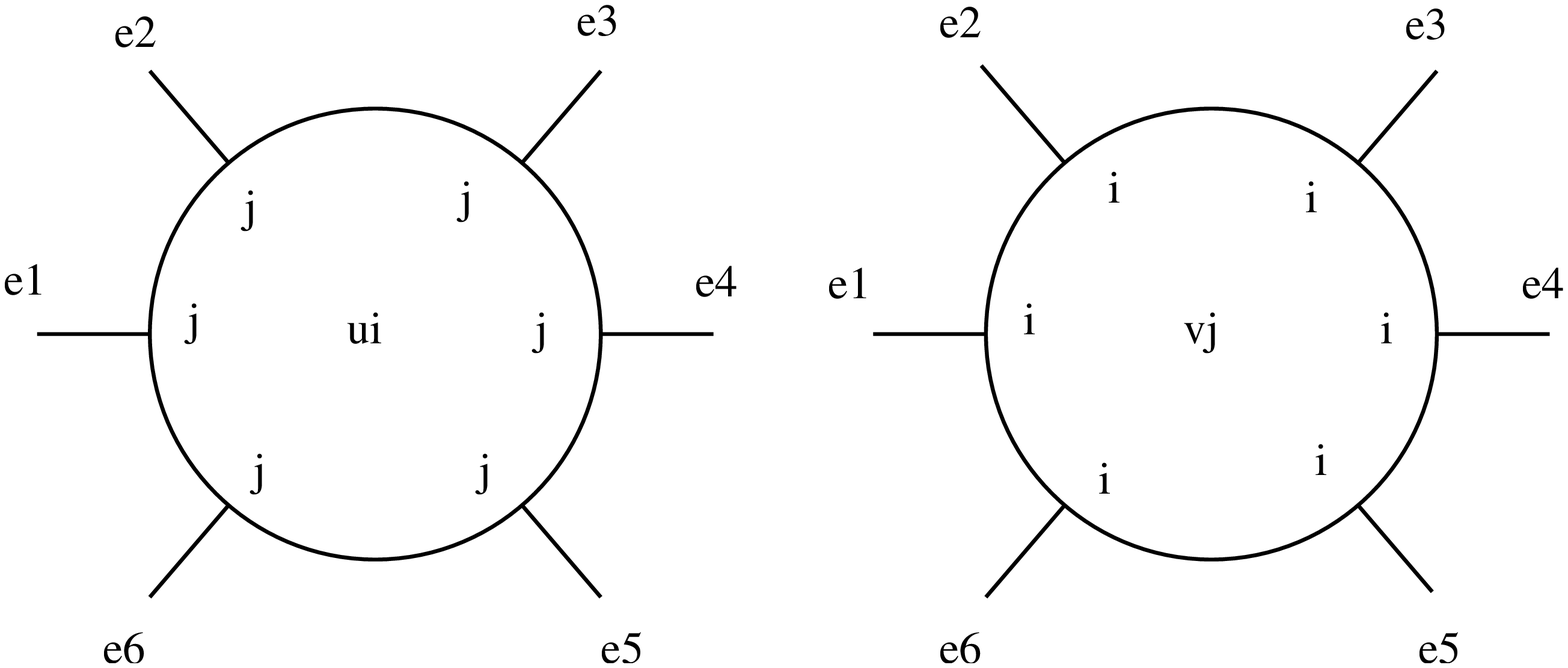}{}{n14b}
\end{figure}

\begin{lem}\label{st3}
The cases $\Delta=6$ and $s,t\geq 3$ do not occur.
\end{lem}

\bpf
Assume $s,t\geq 3$, so that $\Delta=6$ and, in $\bgs,\bgt$,
$\deg\equiv 6$, each edge of $\bgs,\bgt$ has size $t,s$,
respectively, and hence all faces are triangles by
Lemma~\ref{3v}(a). We will say that a vertex in $\bgs,\bgt$ is of
{\it type $(p,n)$} if it has $p$ positive and $n$ negative local
edges, where $p+n=6$.

Suppose some vertex $u_i$ of $\bgs$ is of type $(p,n)$. By the
parity rule, each vertex $v_j$ of $G_T$ has $p$ negative and $n$
positive local edges with label $i$ at $v_j$; thus, by our remarks
above, any vertex of $\bgt$ has $p$ negative and $n$ positive local
edges and is therefore of type $(n,p)$. By a similar argument, any
vertex of $\bgs$ is of type $(p,n)$. Exchanging the roles of $S$ and
$T$ if necessary, we may assume that $(p,n)$ is one of the pairs
$(6,0),(5,1),(4,2),(3,3)$.

If $(p,n)=(6,0)$ then every edge of $\bgt$ is negative, which is
impossible since all the faces of $\bgt$ are triangles, and not all
edges around a triangle face can be negative. Therefore
$(p,n)=(5,1),(4,2),(3,3)$, and so each graph $\bgs,\bgt$ has at
least one positive edge $\ove_S,\ove_T$, of size $t,s$,
respectively.

Now, by Lemma~\ref{pos}, the edge orbits of $\ove_S,\ove_T$ produce
subgraphs isomorphic to the graph in Fig.~\ref{n03}(b) (thick edges
only), and $s,t$ are even. If $\bgs$ has loop edges then $\bgt$ must
have a negative edge $\ove$ which induces the identity permutation;
as $|\ove|=s$, it follows that every vertex of $\bgs$ has an
incident loop edge, and hence that the subgraph of $\bgs$ generated
by the edge orbits of $\ove_T$ and $\ove$ is a union of components
each isomorphic to the graph of Fig.~\ref{n15}(a). Therefore, since
$\deg\equiv 6$ in $\bgs$, if the graph $\bgs$ has loop edges then it
must be of the type shown in Fig.~\ref{n24}(a), where the thick
edges represent the orbits of $\ove_T$. A similar conclusion holds
for $\bgt$ whenever it has any loop edges.
\begin{figure}
\psfrag{p}{$+$}
\psfrag{m}{$-$}
\psfrag{(a)}{$(a)$}
\psfrag{(b)}{$(b)$}
\psfrag{v}{$v$}
\psfrag{v'}{$v'$}
\psfrag{u}{$u$}
\psfrag{u'}{$u'$}
\Figw{5in}{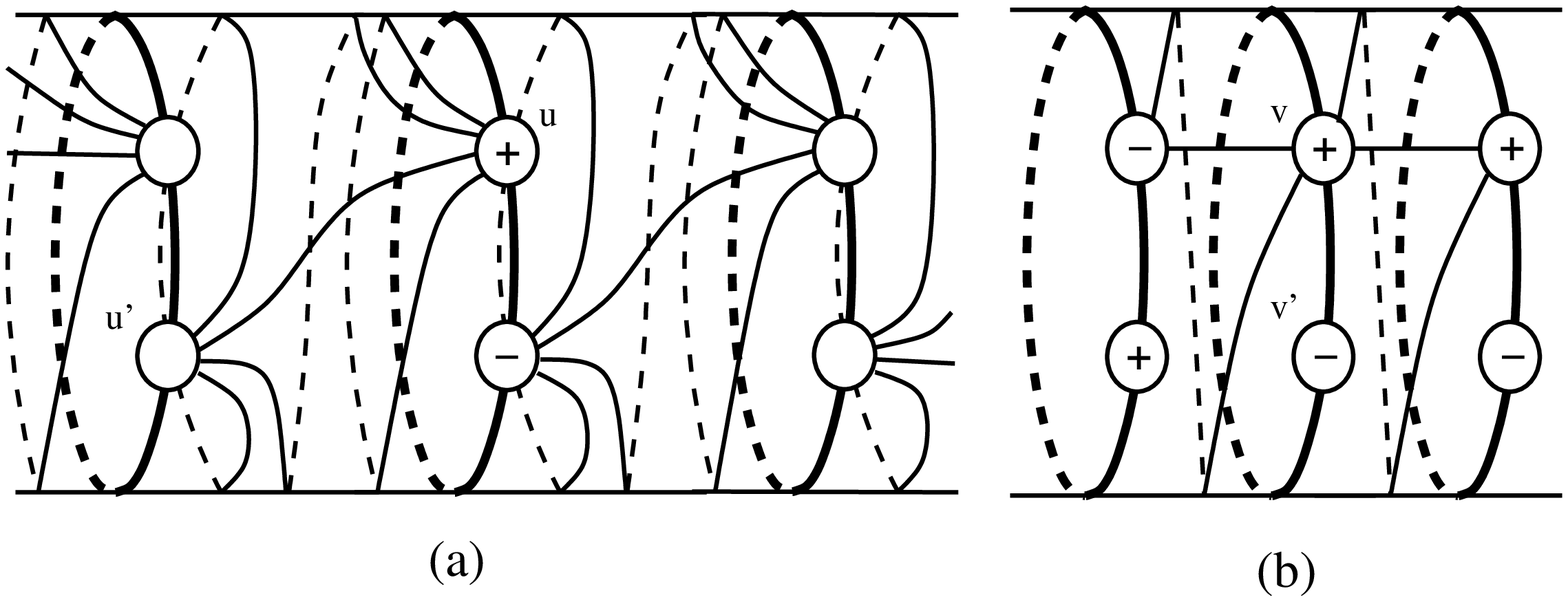}{}{n24}
\end{figure}

If $\bgs$ has no loop edges then we contradict Lemma~\ref{pos},
since each vertex of $\bgs$ has $p\geq 3$ positive local edges. Thus
$\bgs$ has loop edges and so it is a graph of the type shown in
Fig.~\ref{n24}(a). Consider the vertices $u,u'$ of $\bgs$ indicated
in Fig.~\ref{n24}(a), which lie in adjacent edge orbits of $\ove_T$.
If $u$ and $u'$ have opposite parity then $u$ is of type $(2,4)$,
which is not the case. Therefore $u$ and $u'$ have the same parity
and hence are of type $(p,n)=(4,2)$, and the signs of the local
edges as read consecutively around $u$ are of the form $- -++++$. By
the parity rule and the JN1 condition, the signs of the local edges
as read consecutively around each vertex of $\bgt$ are then of the
form $+ + - - - -$.

Let $v,v'$ be the vertices in some edge orbit $c$ of $\ove_S$; then
the two negative edges around, say, the vertex $v$, which are not on
$c$, must both lie on the same side of $c$ (see Fig.~\ref{n24}(b)),
which implies that not both $v,v'$ can have incident loop edges and
hence that $\bgt$ has no loop edges by our preceding arguments. But
then the two positive local edges at $v$ must lie on the other side
of $c$, as shown in Fig.~\ref{n24}(b), so that $\dbgt(v')\leq 4$,
contradicting the fact that $\deg\equiv 6$ in $\bgt$. The lemma
follows.
\epf

\subsection{The cases $s=2, \ t\geq 3$}

By Lemma~\ref{lesst}, $\Delta=6$ and, in $\bgs$, $\deg\equiv 6$ and
all the edges have size $t$; also, recall that any negative edge of
$\bgt$ has size at most $s+1=3$. In these cases $\bgs$ is
combinatorially isomorphic to the graph shown in Fig.~\ref{n13} (cf
\cite[Lemma 5.2]{gordon5}),
with vertices $u_1,u_2$ and edges labeled $\ove_i$, $1\leq i\leq 6$,
and $|\ove_i|=t$.

As $\ove_1,\ove_2$ are positive loop edges in $\bgs$, it follows
from Lemma~\ref{pos} that $t$ is even, so $t\geq 4$. We set $\ve=+1$
if $u_1,u_2$ have the same parity (ie, if $S$ is polarized), and
$\ve=-1$ otherwise (if $S$ is neutral). Then, for $3\leq i\leq 6$,
the edges $\ove_i$ have the same sign as $\ve$.
\begin{figure}
\psfrag{e1}{$\ove_1$}
\psfrag{e2}{$\ove_2$}
\psfrag{e3}{$\ove_3$}
\psfrag{e4}{$\ove_4$}
\psfrag{e5}{$\ove_5$}
\psfrag{e6}{$\ove_6$}
\psfrag{u1}{$u_1$}
\psfrag{u2}{$u_2$}
\psfrag{t}{$t$}
\psfrag{1}{$1$}
\Figw{3.5in}{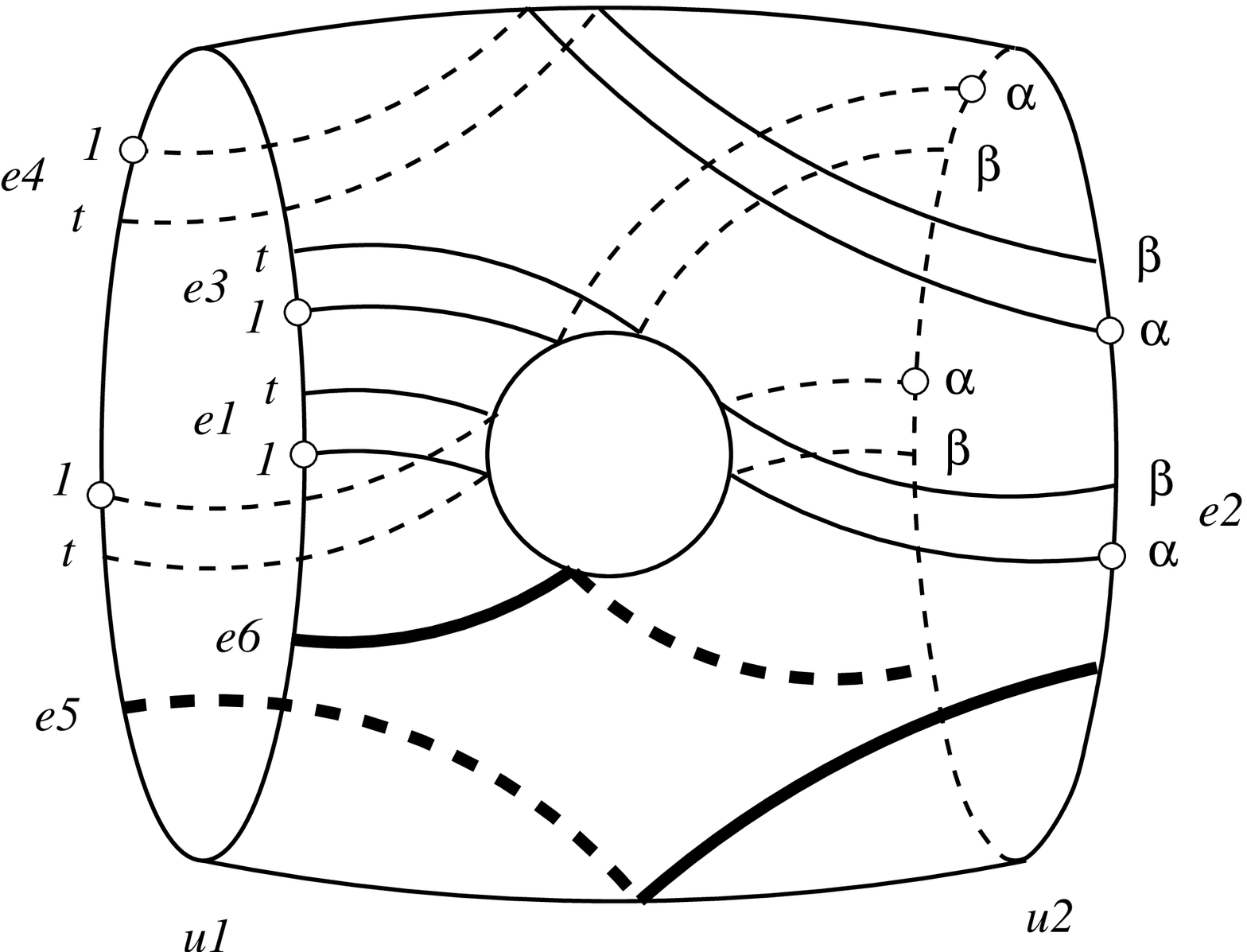}{}{n13}
\end{figure}

Let $\sigma_1,\sigma_2$ be the permutations induced by the edges
$\ove_1,\ove_2$, respectively; notice that
$\ove_3,\ove_4,\ove_5,\ove_6$ all induce the same permutation
$\sigma$. Using the generic labeling scheme of Fig.~\ref{n13} (for
some integers $1\leq \alpha,\beta\leq t$) we can see that
$\sigma_1(x)\equiv 1-x$, $\sigma_2(x)\equiv\alpha+\beta-x$, and
$\sigma(x)\equiv\alpha+\ve-\ve\cdot x \mod t$ for all $1\leq x\leq
t$. As $\alpha+\ve=\sigma(t)=\beta$, we can write $\sigma_2(x)\equiv
2\alpha+\ve-x\mod t$.

The JN1 condition now implies that the local edges around $u_i$ and
$v_j$ for $i=1,2$ and $j=1,\dots, t$ are distributed as shown in
Fig.~\ref{n14} (up to reflections of the vertices).
\begin{figure}
\psfrag{ei}{$e_{\ell}$}
\psfrag{e3}{$e_3$}
\psfrag{e4}{$e_4$}
\psfrag{e5}{$e_5$}
\psfrag{e6}{$e_6$}
\psfrag{ui}{$u_i$}
\psfrag{vj}{$v_j$}
\psfrag{i}{$i$}
\psfrag{j}{$j$}
\Figw{3in}{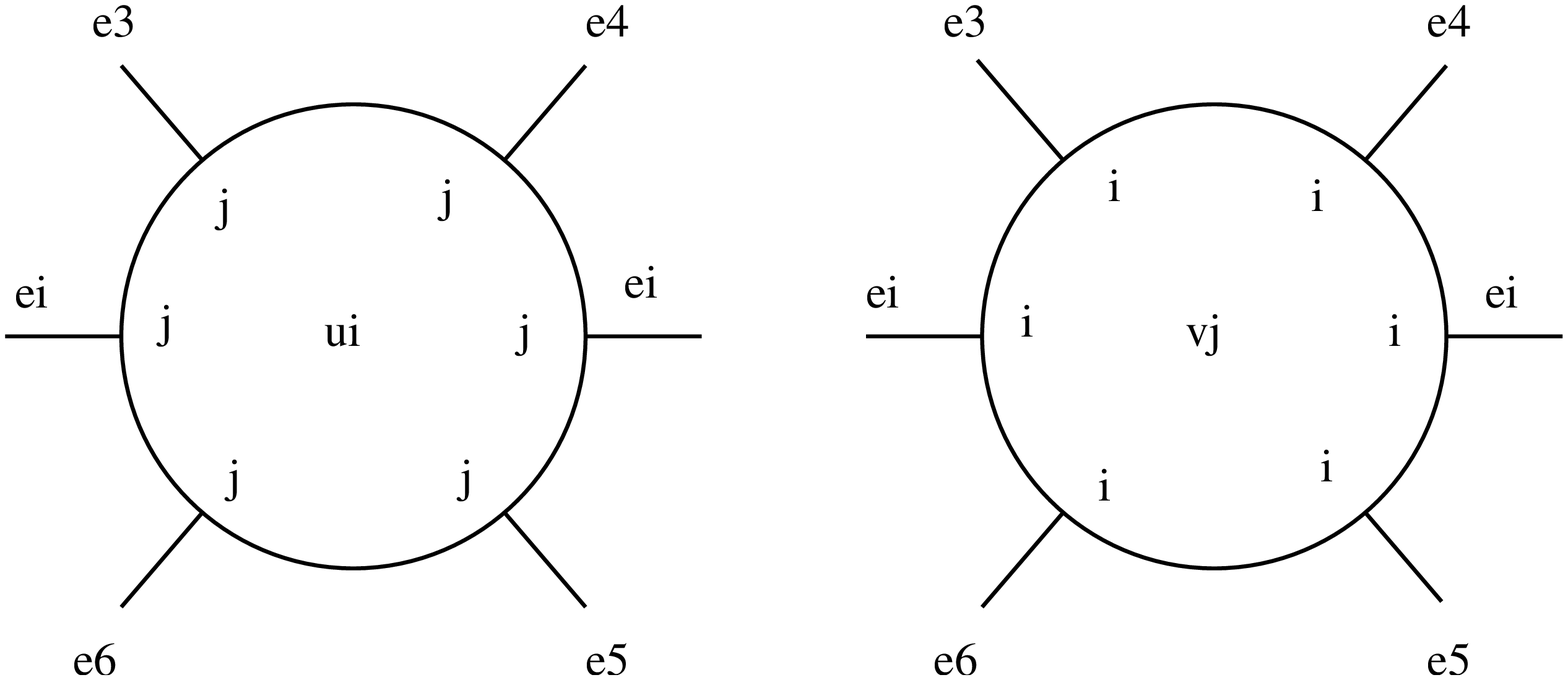}{}{n14}
\end{figure}
In the figure, the edges labeled $e_{\ell}$, $\ell=1,2$, or
$e_3,\dots,e_6$ are edges in the corresponding collections $\ove_k$
of $\bgs$, and represent the same edges in both graphs.

Notice also that some of the local edges around $v_j$ in
Fig.~\ref{n14} may come from distinct parallel edges of $G_T$, since
$\deg\equiv 6$ need not hold in $\bgt$.

\begin{lem}\label{s2a}
The cases $s=2$ and $t\geq 3$ do not occur.
\end{lem}

\bpf
Recall that $t$ is even, so $t\geq 4$. Let $\Gamma$ be the subgraph
of $\bgt$ generated by the edge orbits of $\ove_1$ and the
$\ove_i$'s, $3\leq i\leq 6$.

Observe that $\sigma=\sigma_1
\text{ \ iff \ }\sigma=\sigma_2 \text{  \ iff
\ }\ve=+1 \text{ and }
\alpha\equiv 0\mod t.$ So,
if $\sigma=\sigma_1=\sigma_2$ then $\bgt$ is isomorphic to the
subgraph of $G_T$ generated by the edge orbits of $\ove_1$, and
hence each edge in $\bgt$ is negative of size $6>s+2=4$,
contradicting the fact that any negative edge in $\bgt$ can have
size at most $s+1=3$. Therefore $\sigma\neq\sigma_1,\sigma_2$.

If $\sigma=\id$ then $\Gamma$ is a union of components each
isomorphic to the graph shown in Fig.~\ref{n15}(a), which violates
the JN1 condition. Therefore $\sigma\neq \id$.

Consider any two consecutive cycle edge orbits $\gamma,\gamma'$ of
$\ove_1$ in $G_T$, with opposite parity pairs of vertices $v,v'$ and
$w,w'$, respectively, and denote by $A$ the annular region of $T$
they cobound (see Fig.~\ref{n15}(b)). Let $E$ be the collection of
edges from $\ove_i$, $3\leq i\leq 6$, that lie in $A$. Since
$\sigma\neq\sigma_1,\id$, none of the edges in $E$ are loop edges
nor parallel to the edges in $\gamma,\gamma'$, hence any such edge
has one endpoint on a vertex of $\gamma$ and the other on a vertex
of $\gamma'$.
\begin{figure}
\psfrag{g}{$\gamma$}
\psfrag{g'}{$\gamma'$}
\psfrag{A}{$A$}
\psfrag{(a)}{$(a)$}
\psfrag{(b)}{$(b)$}
\psfrag{(c)}{$(c)$}
\psfrag{e1}{$\ove_1$}
\psfrag{e3}{$a$}
\psfrag{e4}{$b$}
\psfrag{w}{$w'$}
\psfrag{w'}{$w$}
\psfrag{v}{$v$}
\psfrag{v'}{$v'$}
\psfrag{e3'}{$\ove_i$}
\Figw{4.5in}{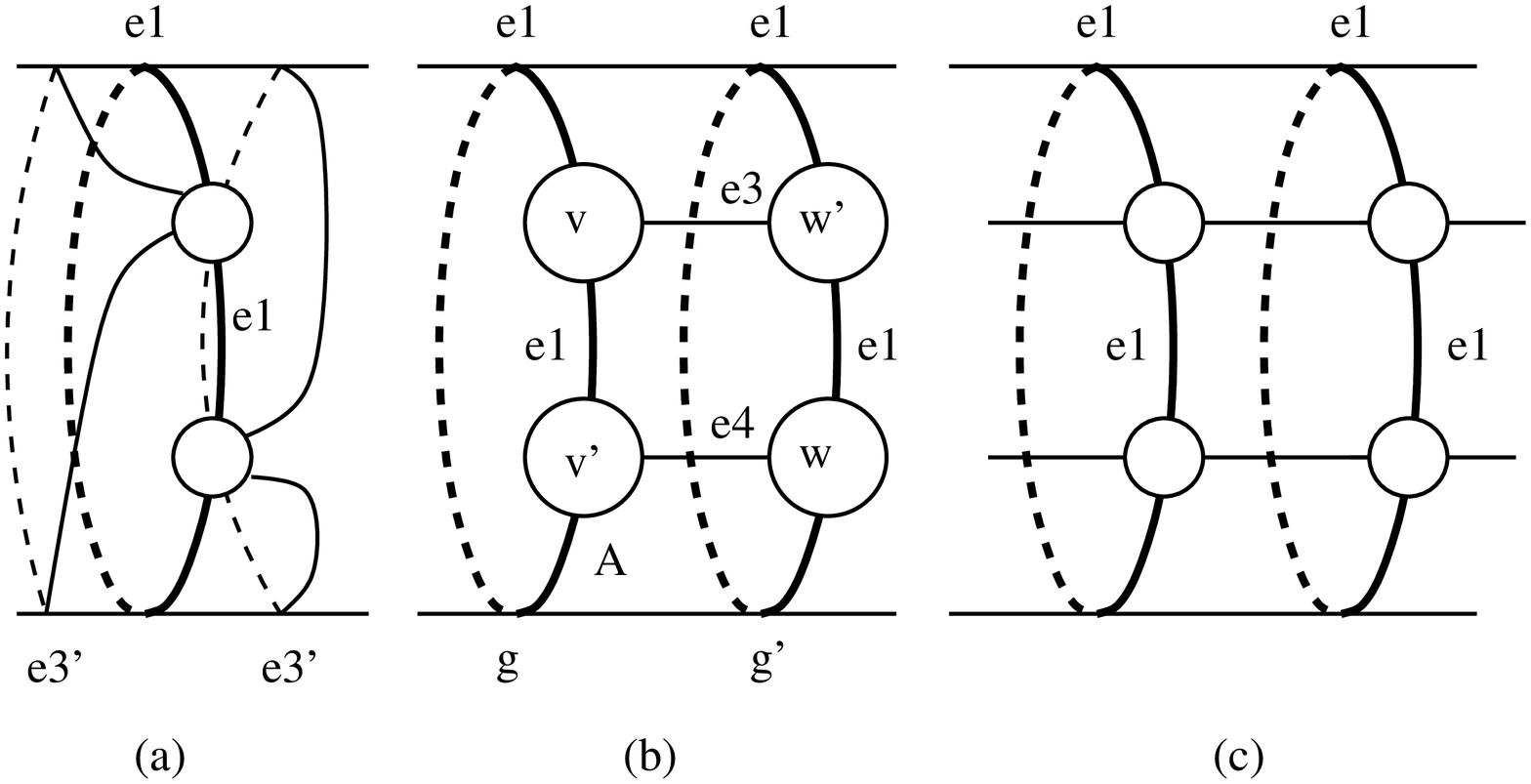}{}{n15}
\end{figure}

Now, by the JN1 condition, each vertex $v,v',w,w'$ has local edges
arising from the edges in $E$. Suppose $a$ is edge of $E$, say with
one endpoint on $v$ and the other on $w$. If $b$ is any edge of $E$
with one endpoint on $v'$ then, by the parity rule, the other
endpoint of $b$ must lie on $w'$ (see Fig.~\ref{n15}(b)). It follows
that the subgraph $\Gamma$ of $\bgt$ is isomorphic to the graph
shown Fig.~\ref{n15}(c), where necessarily each horizontal edge has
size $4=s+2$ and, by Lemma~\ref{mainref}(b), consists of one edge
from each collection $\ove_i$, $3\leq i\leq 6$. Since any negative
edge of $\bgt$ can have size at most $s+1=3$, the horizontal edges
of $\bgt$ must be positive, hence $\ve=-1$ by the parity rule and so
both $S$ and $T$ are neutral.

Thus, any positive edge of $\bgt$ has size 4 and hence its edges
cobound three S-cycle faces in $G_T$, the outermost two of which
locally lie on the same side of $S$ and, by Lemma~\ref{mainref}(b),
have non parallel boundary circles in the surface $S\cup I_{1,2}$ or
$S\cup I_{2,1}$, as the case may be. Therefore, by Lemma~\ref{kb2},
$S$ is generated by an essential once punctured Klein bottle $P$.
\begin{figure}
\psfrag{a}{$\overline{a}$}
\psfrag{b}{$\overline{b}$}
\psfrag{c}{$\overline{c}$}
\psfrag{(a)}{$(a)$}
\psfrag{(b)}{$(b)$}
\Figw{4.3in}{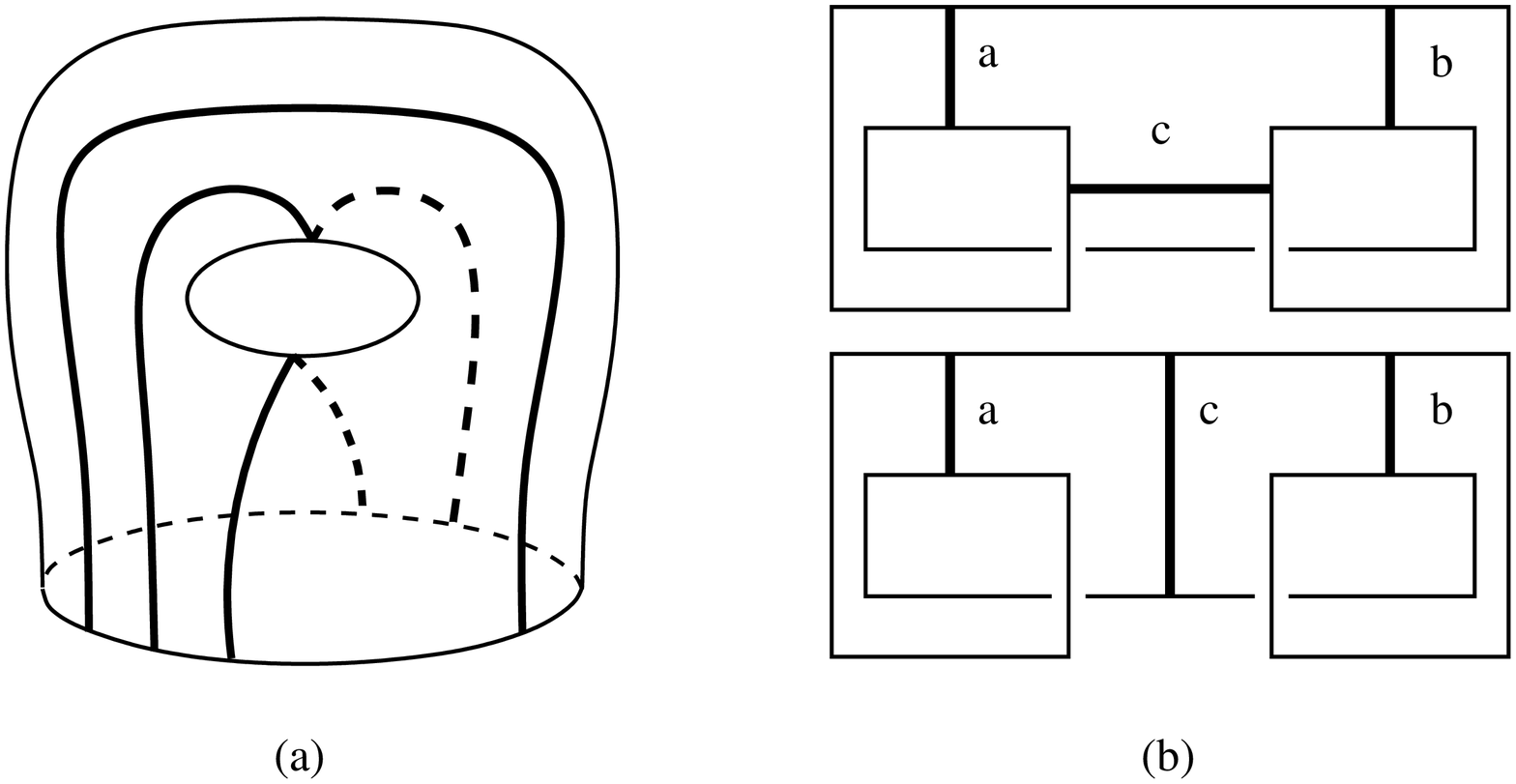}{}{n18}
\end{figure}

Isotope $P$ in $M$ so as to intersect $T$ transversely in essential
graphs. We may assume that $S$ is isotoped accordingly, so that the
new graph $S\cap T\subset T$ is essential and coincides with the
frontier of a small regular neighborhood of the essential graph
$P\cap T\subset T$; in particular, the arguments above apply to the
new graphs $S\cap T\subset T$ and $S\cap T\subset S$. As $P$ has at
most two isotopy classes of negative edges and at most one isotopy
class of positive edges (cf \cite[Lemma 11]{tera1} or
\cite[\S 2]{valdez6}), the reduced graph $\bgpp$ of $G_P=P\cap
T\subset P$ must be isomorphic to one of the graphs in
Fig.~\ref{n18}(b), where the edges $\overline{a},\overline{b}$ are
negative and $\overline{c}$ is positive. Since all the edges of the
reduced graph $S\cap T\subset S$ have size $t$ by Lemma~\ref{lesst},
the edges of $\bgpp$ all have size $t$ too (cf \S\,\ref{ktori}). But
then it is not hard to see that any negative edge in $\bgpp$ induces
the identity permutation, which implies that any negative edge of
the reduced graph of $S\cap T\subset S$ also induces the identity
permutation, contradicting our arguments above on the permutation
$\sigma$. The lemma follows.
\epf

\subsection{The cases $s=1, \ t\geq 3$}
\begin{lem}\label{s1}
The cases $s=1,\ t\geq 3$ do not occur.
\end{lem}

\bpf
If $s=1$ then, as $\deg\equiv 6$ in $\bgs$ by Lemma~\ref{lesst},
$\bgs$ is isomorphic to the graph in Fig.~\ref{n18}(a), and hence
all its edges induce the same permutation $x\mapsto 1-x\mod t$.
Thus, if $\ove$ is any edge of $\bgs$, then $\bgt$ is isomorphic to
the subgraph of $G_T$ generated by the cycle edge orbits of $\ove$,
and so in $\bgt$ each edge is negative of size $3=s+2$,
contradicting the fact that any negative edge in $\bgt$ can have
size at most $s+1=2$.
\epf

\subsection{Proof of Theorem~\ref{main}}

Suppose  $(M,T_0)$ is not cabled and $(F_1,\partial F_1)$,
$(F_2,\partial F_2)\subset (M,T_0)$ are $\mc{K}$-incompressible tori
with boundary slopes at distance $\Delta\geq 6$. We set
$\{S,T\}=\{F_1,F_2\}$, with $s=|\partial S|$ and $t=|\partial T|$.
Then $1\leq s,t\leq 2$ by Lemmas~\ref{st3},
\ref{s2a}, and \ref{s1}, so it only remains to check that
$\Delta\leq 8$. The case $s=t=1$ is impossible by the parity rule,
and there are three more cases to consider.

\setcounter{case}{0}
\begin{case}
$s=t=2$ and $S$ is polarized.
\end{case}

Then $T$ is neutral and all the edges of $G_T$ are negative. Hence
$\bgt$ has at most 4 edges, each of size at most $s+1=3$, and
$\deg\leq 4$ in $\bgt$, so the degree of $v_1$ in $G_T$ satisfies
the relations $2\Delta=s\cdot\Delta=\dgt(v_1)\leq 4\cdot 3=12$. Thus
$\Delta=6$ (and $\deg\equiv 4$ in $\bgt$, with each edge of $\bgt$
of size $s+1=3$).

\begin{case}
$s=t=2$ and both $S$ and $T$ are neutral.
\end{case}

Then, in either graph $\bgs,\bgt$, any vertex has at most 4 negative
edges and either 0 or 2 positive local edges (see Fig.~\ref{n13}),
hence by Lemma~\ref{mainref}(b) any local positive edge has size at
most 4, while any negative edge has size at most 2. Therefore, if
$p,n$ are the number of positive and negative local edges of $\bgs$
at $u_1$, respectively, then $p\leq 2$ and $n\leq 4$, so the degree
of $u_1$ in $G_S$ satisfies the relations
$2\Delta=s\cdot\Delta=\dgs(u_1)\leq p\cdot 4+n\cdot 2\leq 16,$ and
so $\Delta\leq 8$.

\begin{case}
$s=1$ and $t=2$.
\end{case}

By the parity rule, since $S$ is polarized then $T$ is neutral and
all edges of $G_T$ negative, so $\bgt$ has at most 4 edges, each of
size at most $s+1=2$. Hence $\deg\leq 4$ in $\bgt$, and so the
degree of $v_1$ in $G_T$ satisfies the relations
$\Delta=s\cdot\Delta=\dgt(v_1)\leq 4\cdot 2\leq 8$.
\hfill\qed


\providecommand{\bysame}{\leavevmode\hbox to3em{\hrulefill}\thinspace}
\providecommand{\MR}{\relax\ifhmode\unskip\space\fi MR }
\providecommand{\MRhref}[2]{%
  \href{http://www.ams.org/mathscinet-getitem?mr=#1}{#2}
}
\providecommand{\href}[2]{#2}


\begin{thebibliography}{10}

\bibitem{cgls}
M.~Culler, C.~McA. Gordon, J.~Luecke, and P.~Shalen, \emph{Dehn
surgery on
  knots}, Ann. of Math. (2) \textbf{125} (1987), no.~2, 237--300.
  \MR{88a:57026}

\bibitem{gordon5}
C.~McA. Gordon, \emph{Boundary slopes of punctured tori in
$3$-manifolds},
  Trans. Amer. Math. Soc. \textbf{350} (1998), no.~5, 1713--1790.
  \MR{98h:57032}

\bibitem{gordonlith}
C.~McA. Gordon and R.~A. Litherland, \emph{Incompressible planar
surfaces in
  {$3$}-manifolds}, Topology Appl. \textbf{18} (1984), no.~2-3, 121--144.
  \MR{86e:57013}

\bibitem{gordonlu4}
C.~McA. Gordon and J.~Luecke, \emph{Toroidal and boundary-reducing
{D}ehn
  fillings}, Topology Appl. \textbf{93} (1999), no.~1, 77--90. \MR{2000b:57030}

\bibitem{tera1}
Kazuhiro Ichihara, Masahiro Ohtouge, and Masakazu Teragaito,
\emph{Boundary
  slopes of non-orientable {S}eifert surfaces for knots}, Topology Appl.
  \textbf{122} (2002), no.~3, 467--478. \MR{1 911 694}

\bibitem{lee1}
Sangyop Lee, \emph{Exceptional Dehn fillings on hyperbolic
3-manifolds with at
  least two boundary components}, Prepint, 2005.

\bibitem{tera10}
Sangyop Lee, Seungsang Oh, and Masakazu Teragaito, \emph{{Reducing
Dehn
  fillings and small surfaces}}, Proc. London Math. Soc. \textbf{92} (2006),
  203--223.

\bibitem{mati1}
Daniel Matignon and Nabil Sayari, \emph{Klein slopes on hyperbolic
3-manfolds},
  Preprint, 2005.

\bibitem{valdez6}
Enrique Ram{\'{\i}}rez-Losada and Luis~G.\ Valdez-S{\'a}nchez,
  \emph{Once-punctured {K}lein bottles in knot complements}, Topology Appl.
  \textbf{146/147} (2005), 159--188. \MR{MR2107143 (2005i:57009)}

\bibitem{tera6}
Masakazu Teragaito, \emph{Distance between toroidal surgeries on
hyperbolic
  knots in the 3-sphere}, Trans. Amer. Math. Soc. \textbf{358} (2006), no.~3,
  1051--1075.

\end{thebibliography}
\end{document}